\documentclass{amsart}
\usepackage{amssymb, amsmath, amsfonts, amscd}

\input xypic

\theoremstyle{plain}
\newtheorem{theorem}{Theorem}[section]
\newtheorem{corollary}[theorem]{Corollary}
\newtheorem{lemma}[theorem]{Lemma}
\newtheorem{proposition}[theorem]{Proposition}
\newtheorem{conjecture}{Conjecture}
\newtheorem{example}[theorem]{Example}

\theoremstyle{definition}
\newtheorem{definition}[theorem]{Definition}

\theoremstyle{remark}

\numberwithin{equation}{theorem}

\newcommand{\A}{\mathcal{A}}

\newcommand{\E}{\mathcal{E}}

\renewcommand{\O}{\mathcal{O} }

\newcommand{\Pic}{\operatorname{Pic} }
\newcommand{\CH}{\operatorname{CH} }

\newcommand{\Spec}{\operatorname{Spec} }
\renewcommand{\P}{\operatorname{P} }
\renewcommand{\H}{\operatorname{H} }

\newcommand{\K}{\operatorname{K} }
\newcommand{\Z}{\mathbb{Z} }

\renewcommand{\K}{\operatorname{K} }

\renewcommand{\L}{\operatorname{L} }

\newcommand{\Proj}{\operatorname{Proj}}

\newcommand{\C}{\mathbb{C} }
\newcommand{\Q}{\mathbb{Q} }
\newcommand{\Soule}{Soul\'{e} }
\renewcommand{\P}{\mathbb{P}}
\newcommand{\G}{\mathbb{G}}
\renewcommand{\A}{\mathbb{A}}
\newcommand{\Fl}{\mathbb{F}} 
\newcommand{\Log}{\operatorname{Log}}
\renewcommand{\ln}{\operatorname{ln}}

\newcommand{\SL}{\operatorname{SL}}

\begin{document}

\title{Conjectures on L-functions for flag bundles on Dedekind domains}

\author{Helge \"{O}ystein Maakestad}

\email{\text{h\_maakestad@hotmail.com}}

\thanks{}


\date{June 2020}

\begin{abstract}  The aim of this paper is to give evidence for the Beilinson-\Soule vanishing conjecture and \Soule conjecture on L-functions for a class of bundles over Dedekind domains.
Let $\O_K$ be the ring of integers in an algebraic number field $K$ with $S:=\Spec(\O_K)$. Let $T_0,\ldots, T_n$ be regular schemes of finite type over $S$, and let $X$ be a scheme of finite type over $T_n$
with a stratification (a generalized cellular decomposition) of closed subschemes
\[ \emptyset=X_{-1} \subseteq X_0 \subseteq \cdots \subseteq X_{n-1} \subseteq X_n:=X \]
with $X_i-X_{i-1}=E_i$, where $E_i$ is a vector bundle of rank $d_i$ on $T_i$.  We prove that if the Beilinson-\Soule vanishing conjecture and the \Soule conjecture on L-functions holds for $T_i$, it follows the same conjectures hold for $X$.  We moreover prove the Beilinson-\Soule vanishing conjecture and the \Soule conjecture on L-functions for any partial flag bundle $\mathbb{F}(N,\E)$ where $\E$ is a coherent $\O_S$-module.  We also prove the conjectures for any finite rank vector bundle and affine fibration on $S$.
We reduce the study of the Beilinson-\Soule vanishing conjecture and the \Soule conjecture on L-functions to the study of affine regular schemes of finite type over $\Z$. Hence we get an approach to the Birch and Swinnerton-Dyer conjecture for abelian schemes using affine regular schemes of finite type over $\Z$. For a partial flag bundle $\Fl(N,E)$ over $\O_K$ we give an explicit formula for the L-funtion in terms of the L-function of $\O_K$. Hence the Bloch-Kato conjecture on the Tamagawa number of $\Fl(N,E)$ is reduced to the study of the same conjecture for $\O_K$.
\end{abstract}

\maketitle

\tableofcontents

\section{Introduction} If $\O_K$ is the ring of integers in an algebraic number field $K$, it follows the rank of the m'th K-group $\K'_m(\O_K)$ and the rank of the weight space $\K'_m(\O_K)_{(i)}$ 
is well known for all integers $m \geq 0, i\geq 1$  (see \cite{Borel} and \cite{Kahn}). We may define the L-function $\L(\O_K,s)$ of $\O_K$ and the K-theoretic j'th Euler characteristic
\begin{align}
&\label{eulerint}  \chi(\O_K, j):= \sum_{m\geq 0}(-1)^{m+1}dim_{\Q}(\K'_m(\O_K)_{(j)}) .
\end{align}
The function $\L(\O_K,s)$ is the well known Dedekind L-function of the number  field $K$.
Borel proved in \cite{Borel} that 
\begin{align}
&\label{eulerintro} \text{$\chi(\O_K,j)$ is an integer for any number field $K$ and any integer $j$.}
\end{align}
He also proved the relationship
\begin{align}
 &\label{orderintro} \chi(\O_K, j)=ord_{s=j}(\L(\O_K, s)) 
\end{align}
between the Euler characteristic and L-function of $\O_K$. 
In the litterature Conjecture \ref{eulerintro} is referred to as the Beilinson-\Soule vanishing conjecture, and Conjecture \ref{orderintro} is referred to as the \Soule conjecture. In the paper \cite{Soule} Conjectures \ref{eulerintro} and \ref{orderintro} are formulated for any quasi projective scheme of finite type over $\Z$. The aim of this paper is to prove the Beilinson-\Soule vanishing conjecture and the \Soule conjecture  for 
a class of schemes of finite type over the ring $\O_K$ called partial flag bundles.

Let $k$ be a field and $W$ an $n$-dimensional vector space over $k$. Let $N:=\{n_1,..,n_l\}$ be a sequence of positive integers with $\sum_i n_i=n$ and let $\Fl(N,W)$ be the flag variety of flags of type $N$ in 
$W$. It follows the set of $k$-rational points of $\Fl(N,W)$ are in one-to-one correspondence with the set of flags $\{W_i\}$ of type $N$ in $W$. A flag of type $N$ in $W$ is a sequence of $k$-vector spaces 
\[ W_1\subseteq W_2 \subseteq \cdots \subseteq W_{l-1} \subseteq W \]
with $dim_k(W_i)=n_1+\cdots +n_i$. If $l=2$ and $n_1<n$, it follows the flag variety $\Fl(N,W)$ is the grassmannian variety $\G(n_1,W)$ of $n_1$-dimensional sub spaces of $W$.
A partial flag bundle is a relative version of $\Fl(N,W)$. Let $S$ be a scheme and let $\E$ be a coherent $\O_S$-module. The flag bundle $\Fl(N,\E)$ is a scheme equipped with a  surjective morphism of schemes
\[ \pi: \Fl(N,\E) \rightarrow S \]
such that the fiber $\pi^{-1}(s)$ at any point $s \in S$ is isomorphic to the flag variety $\Fl(N,\E(s))$ of flags of $\kappa(s)$-vector spaces of type $N$ in the fiber $\E(s)$ of $\E$ at $s$. If $\E$ is locally trivial it follows the map $\pi: \Fl(N,E)\rightarrow S$ 
is locally trivial in the Zariski  topology.

Let $T$ be a regular scheme of finite type over $\O_K$ and $X$ be a scheme of finite type over $T$ with a cellular decomposition $X_i\subseteq X$ of closed subschemes, such that $X_i-X_{i-1}$ is a finite disjoint union of affine space $\A^i_T$ over $T$. In Theorem \ref{specialmain} we prove that if the Beilinson-\Soule vanishing conjecture and the \Soule conjecture on L-functions holds for $T$, it follows the same conjectures hold for $X$.  In particular it follows there is an equality of integers
\[ \chi(X,j)=ord_{s=j}(\L(X,s))\]
(see Theorem  \ref{specialmain}). We prove a similar result for any scheme $X$ equipped with a generalized cellular decomposition in Lemma \ref{lemmagen}.
We moreover prove the Beilinson-\Soule vanishing conjecture and the \Soule conjecture on L-functions for any partial flag bundle $\mathbb{F}(N,\E)$ on $\Spec(\O_K)$ where $\E$ is a coherent $\O_S$-module (see Corollary \ref{s13corollary}). This gives an infinite number of non-trivial examples of partial flag bundles $\Fl(N,\E)$ where Conjecture \ref{sou1} and \ref{sou3} hold (see Example \ref{nontrivial}).
If $A$ is an abelian scheme over $\O_K$ it follows the \Soule Conjecture \ref{sou3} for $A$ is one way to formulate a version of the \emph{Birch and Swinnerton-Dyer conjecture} for $A$ using algebraic K-theory.

We reduce the study of the Beilinson-\Soule vanishing conjecture and the \Soule conjecture on L-functions to the study of affine regular schemes of finite type over $\Z$.
Hence we get an approach to the Birch and Swinnerton-Dyer conjecture for abelian schemes using affine regular schemes of finite type over $\Z$.

We use the projective bundle formula for algebraic K-theory and an elementary construction of eigenvectors for the Adams operator to calculate the weight space 
\[ \K'_m(\P(E^*))_{(i)} \]
for any pair of integers $m\geq 0, i\geq 1$, any finite rank projective $\O_K$-module $E$ for any algebraic number field $K$ (see Theorem \ref{mainthm}).
This illustrates the possiblility to do explicit computations for the K-theory of projective bundles and more general flag bundles.

In Chapter two we introduce some notation and state known results on algebraic K'-theory and Adams operations for rings of integers in algebraic number fields. In Chapter three we state some conjectures on the algebraic K'-theory
$\K'_m(X)_{\Q}^{(i)}$ of a quasi projective scheme $X$ over the integers $\Z$ and the Hasse-Weil L-function $\L(X,s)$ of $X$. In Chapter four we prove the mentioned conjectures for any scheme equipped with a generalized cellular decomposition.
As a particular case we prove the conjectures for any complete flag bundle over the ring of integers $\O_K$ in any number field $K$. In Appendix A we give an elementary proof of an explicit formula for the i'th Adams weight space $\K_m'(\P(E^*))_{\Q}^{(i)}$
of the algebraic K'-theory of any projective bundle $\P(E^*)$ on the ring of integers $\O_K$ in any number field $K$. This gives an explicit formula for the i'th Euler characteristic $\chi(\P(E^*),i)$ for any projective bundle $\P(E^*)$ over $\O_K$.
In Appendix B we prove some elementary properties of formal power series that are used in the calculations in Appendix A.


\section{Algebraic K-theory and Adams operations}

Let $\O_K$ be the ring of integers in an algebraic number field $K$ and let $S=\Spec(\O_K)$. Let $X$ be a scheme of finite type over $S$.

In this section we introduce some notation from \Soule's original paper \cite{Soule}: Let $M(X)$ denote the category of coherent $\O_X$-modules and let $BQM(X)$ denote the simplicial classifying set of $M(X)$. Let $BQP(X)$ denote the simplicial classifying set of $P(X)$, where $P(X)$ is the category of locally trivial finite rank $\O_X$-modules. By definition

\begin{align}
\K’_m(X) &\label{kp1}:= \pi_{m+1}(BQM(X)) \\
\K_m(X) &\label{kp2}:= \pi_{m+1}(BQP(X))
\end{align}

where m is an arbitrary integer. If $X$ is a regular scheme it follows $\K’_m(X)=\K_m(X)$ and $\K_m(X)=0$ for $m<0$. Assume $X$ is a scheme of finite type over $\Z$ and assume $u:X \rightarrow M$ is a closed immersion into a scheme $M$ where $M$ is a regular scheme of finite type over $\Z$ of dimension $D$. Define
$\K^X_m(M)$ as the homotopy group of the fiber of the canonical map $BQP(M) \rightarrow BQP(M-X)$


\begin{definition}\label{adamsoperation}
If $Y$ is a regular scheme of finite type over $\Z$, there is for every positive integer $k \geq 0$ an action
\[ \psi^k: \K_m(Y) \rightarrow \K_m(Y) \]
with the following properties: If $L$ is the class of a line bundle in $\K_0(Y)$ it follows 
\[ \psi^k(L):=L^k .\]
The map $\psi^k$ is the \emph{$k$'th Adams operator for $\K_m(Y)$}.
\end{definition}

The map $\psi^k$ is functorial in the sense that for any map $p:Y\rightarrow Y'$ of regular schemes $Y,Y'$ of finite type over $\Z$ it follows
\[ \psi^k(p^*x)=p^*(\psi^k(x)) \]
for any element $x\in \K_m(Y')$. The abelian group $\K_*(Y):=\oplus_{m\geq 0}\K_m(Y)$ is a graded commutative ring and the endomorphism
\[ \psi^k: \K_*(Y) \rightarrow \K_*(Y) \]
is a ring homomorphism: $\psi^k(xy)=\psi^k(x)\psi^k(y)$ for any $x\in \K_m(Y), y\in\K_n(Y)$. The operation $\psi^k$ induce canonically a ring homomorphism
\[ \psi^k:\K_*(Y)\otimes \Q \rightarrow \K_*(Y)\otimes \Q \]
(let $\K_m(Y)_{\Q}:=\K_m(Y)\otimes \Q$) and we define 
\[ \K_m(Y)_{\Q}^{(i)} :=\{ x\in \K_m(Y)_{\Q}:\text{ such that $\psi^k(x)=k^ix$.} \}\]

There is a direct sum decomposition
\[ \K_m(Y)_{\Q}\cong \oplus_{i\in \Z}K_m(Y)_{\Q}^{(i)} \]
and the space $\K_m(Y)_{\Q}^{(i)}$ is independent of choice of positive integer $k$.
By definition we let
\begin{align}
&\label{kth2} \K’_m(X)_{(i)} := \K^X_m(M)_{\Q}^{(D-i)}. 
\end{align}






When $X$ is regular we may choose $M=X$. It follows

\[ \K’_m(X)_{(i)} =\K_m(X)_{\Q}^{(D-i)}\]

where $D=dim(X)$. Hence when $X$ is a regular scheme of finite type over $\Z$ we may use the K-theory of the category $P(X)$ of finite rank algebraic vector bundles on $X$ and the 
Adams operations on $\K_m(X)_{\Q}$ to calculate the group $\K’_m(X)_{(i)}$ introduced in \Soule’s paper.

\begin{definition} Let $X$ be a scheme of finite type over $\Z$ and let $i: X\rightarrow M$ be a closed embedding into a regular scheme $M$ of finite type over $\Z$ with $D:=dim(M)$.
Define
\[ \K'_*(X):=\oplus_{m\in \Z}\K_m^X(M) \]
and $\K'_*(X)_{\Q}:=\K'_*(X)\otimes \Q$. Define $\K'_m(X)_{(j)}:=\K_m^X(M)_{\Q}^{(D-j)}$.
The $\Q$-vector space $\K'_m(X)_{(j)}$ is the \emph{weight space of weight $j$}.
\end{definition}

The following result calculates $\K_m(\O_K)_{\Q}$ and $\K_m(\O_K)_{\Q}^{(i)}$ for all $m,i$:

\begin{theorem} \label{boreltheorem} Let $K$ be a number field with ring of integers $\O_K$ and real and complex places $r_1,r_2$. The following holds:
\begin{align}
&\label{borm}\K_m(\O_K)_{\Q}=0\text{ for all $m<0$}\\
&\label{bor0}\K_0(\O_K)_{\Q}=\Q\\
&\label{bor1}\K_m(\O_K)_{\Q}=0\text{ for $m=2i, i\neq 0$}\\
&\label{bor2}\K_m(\O_K)_{\Q}=\Q^{r_1+r_2}\text{ for $m\equiv 1 \text{ mod }4$}\\
&\label{bor3}\K_m(\O_K)_{\Q}=\Q^{r_2}\text{ for $m\equiv 3 \text{ mod }4$}.
\end{align}
Moreover
\begin{align}
&\label{bor4}\K_{2i-1}(\O_K)_{\Q}^{(i)} = \Q^{r_1+r_2}\text{ for $i\equiv 0\text{ mod }2$}\\
&\label{bor5}\K_{2i-1}(\O_K)_{\Q}^{(i)}=\Q^{r_2}\text{ for $i\equiv 1\text{ mod } 2$}.
\end{align}
\end{theorem}

The reader should consult \cite{Kahn} and \cite{soule1} for the history of the calculation of Theorem \ref{boreltheorem}.
The calculation of $\K'_m(\O_K)\otimes \Q$ for $m=0$ follows from the fact the ideal class group of $\O_K$ is finite, a result going back to Minkowski. The case $m=1$ is Dirichlet's unit theorem.
For $m\geq 2$ Theorem \ref{boreltheorem} follows from Borels paper \cite{Borel}.  The formula for the weight space decomposition is proved in \cite{Kahn} and the book also gives references to the papers \cite{soule2} and \cite{soule3}.

\section{Reduction of the Beilinson-\Soule vanishing conjecture and \Soule conjecture on L-functions to the affine regular case}

In this section we give a criteria for Conjecture \ref{sou1} to hold for a scheme $X$ of finite type over $\Z$ in terms of an open cover $U_i$ of $X$. We also study the L-function $\L(X,s)$ and $ord_{s=k}(\L(X,s))$ in 
terms of the cover $U_i$ (see Lemma \ref{eulercover} and \ref{lcover}).  We reduce the study of the Beilinson-\Soule vanishing conjecture and \Soule's conjecture on L-functions to the study of affine regular schemes of finite type over $\Z$.
We also prove the Beilinson-Soule vanishing conjecture and Soule conjecture on L-functions for any affine fibration  $\mathbb{V}(E^*)$ on $\Spec(\O_K)$ where $E$ is a coherent $\O_K$-moduler and $K$ a number field 
(see Theorem \ref{vbtheorem} and \ref{affinefib}).

Let in the following $K$ be an algebraic number field with ring of integers $\O_K$ and let $S:=\Spec(\O_K)$.



\begin{definition}\label{euler} Let $X$ be a quasi projective scheme of finite type over $S$ and let $i\in \Z$. Let 
\[ \chi(X,i):= \sum_{m \in \Z}(-1)^{m+1}dim_{\Q}(\K’_m(X)_{(i)}) \]
be the \emph{Euler characteristic} of $X$ of type $i$.
\end{definition}

Note: The Euler characteristic $\chi(X,i)$ may not be a well defined integer in general since the sum in Definition \label{euler} is infinite.


\begin{definition} \label{lfunction} Let $X$ be a scheme of finite type over $S$.
Let 

\[ \L(X,s) := \prod_{x \in X^{cl}}\frac{1}{1-N(x)^{-s}} \]

be the \emph{Hasse-Weil L-function} of $X$. Here we view $s$ as a complex variable and the infinite 
product is taken over the set of closed points $x$ in $X^{cl}$. By definition $N(x):=\#\kappa(x)$ where $\kappa(x)$ is the residue field of $x$.
\end{definition}

Note: Since $X$ is of finite type over $\Z$ and $x$ is a closed point, it follows $\kappa(x)$ is a finite field. The formal infinite product $\L(X,s)$ is not a well defined function for all complex numbers $s$. There 
are long standing conjectures on convergence properties of $\L(X,s)$ as a complex function in the variable $s$ in the case when $X=\Spec(\O_K)$ with $K$ a number field (see \cite{neukirch}, Chapter VII).

\begin{example} Hasse-Weil L-functions and Weil zeta functions.\end{example}

Let in this example $T:=\Spec(k)$ with $k$ a finite field with $q$ elements and let $\pi: X \rightarrow T$ be a scheme of finite type over $T$. It follows there is for any integer $r \geq 2$ a canonical 
finite extension of fields $k \subseteq k_r$. Let $\overline{k}$ be an algebraic closure of $k$ and let $\overline{X}:=X\times_T \Spec(\overline{k})$. 
Let $N_r:=\#\overline{X}(k_r)$ be the number of points of $\overline{X}$ that are rational over $k_r$ and define

\begin{align}
\label{weil} Z(X,t)&:= \operatorname{exp}(\sum_{r \geq 1}\frac{N_r}{r}t^r).
\end{align}

\begin{definition} 
Let $X$ be a scheme of finite type over $T$. The formal power series $Z(X,t)$ defined in 
\ref{weil} is the  \emph{Weil zeta function} of $X$.
\end{definition}

\begin{lemma} \label{productL}  Let $X$ be a quasi projective scheme of finite type over $T$ and let $q:=\#k$ be the number of elements in $k$. There is a formal equality
\begin{align}
\label{eqprod}  \L(X,s)&= Z(X,q^{-s}) ,
\end{align}
where $\L(X,s)$ is the Hasse-Weil L-function of $X$.
\end{lemma}
\begin{proof} Note: Since $X$ is of finite type over $T$ it follows $X$ is of finite type over $\Spec(\Z)$ and we may defined the L-function $\L(X,s)$.
See \cite{hartshorne}, Exercise 5.4 in Appendix C for a proof of the equality in \ref{eqprod}.
\end{proof}

Note: Since $\L(X,s)$ is a formal infinite product involving the complex variable $s$ and $Z(X,t)$ is a formal power series living in $\Q[[t]]$, one has to give precise meaning to the equality in Lemma \ref{productL} but this is not needed in this paper.

\begin{example} The Dedekind L-function.\end{example}
If $K$ is an algebraic number field with ring of integers $\O_K$ and $S:=\Spec(\O_K)$, it follows $\L(S,s)$ is the \emph{Dedekind L-function} of $K$.
In particular $\L(\Spec(\Z),s)$ is the \emph{Riemann zeta function}.

In \Soule’s paper \cite{Soule} the following conjecture is stated:

\begin{conjecture}\label{souleconjecture} (Conjecture 2.2 in \cite{Soule}) Let $X$ be a quasi projective scheme of finite type over $\Z$ and let $i \in \Z$ be an integer.
\begin{align} 
&\label{sou1}\text{For fixed integer $i$ the group }\K’_m(X)_{(i)}\text{ is zero for almost all integers $m$.}\\
&\label{sou2}dim_{\Q}(\K’_m(X)_{(i)}) \text{ is finite for all $m, i$.}\\
&\label{sou3}\chi(X,i)=ord_{s=i}(\L(X,s))\text{ for all $i\in \Z$} 
\end{align}
\end{conjecture}

Note: The Conjecture \ref{sou3} is mentioned in Wiles' CLAY Math description of the \emph{Birch and Swinnerton-Dyer conejcture} (one of the \emph{Millenium Problems}, see \cite{wiles}). 
In \cite{wiles} Conjecture \ref{sou3} is referred to as due to Tate, Lichtenbaum, Deligne, Bloch, Beilinson and others. 
Conjecture \ref{sou1} is sometimes referred to as the \emph{Beilinson-\Soule vanishing conjecture}. 
If $E$ is a relative elliptic curve over $\O_K$, it follows Conjecture \ref{sou3} is a version of the Birch and Swinnerton-Dyer conjecture for $E$ using K-theory.  
The version given in \cite{wiles} is formulated for an elliptic curve $E$ over $\Q$ and the group of rational points $E(\Q)$ of $E$. 
There is an embedding $E(\Q)\subseteq \operatorname{Pic}(E)$ and $\K_0(E)=\operatorname{Pic}(E)\oplus \Z$, hence the conjecture in \cite{wiles}
is similar to Conjecture \ref{sou3}. Hence we may view the conjecture mentioned in \cite{wiles}  as a special case of Conjecture \ref{sou3}.
Note morover that if $X_{red}$ is the reduced scheme of $X$ it follows $\L(X,s)=\L(X_{red},s)$ and $\K'_m(X)=\K'_m(X_{red})$, hence Conjecture \ref{souleconjecture} holds for $X$ if and only if it holds for 
$X_{red}$.

\begin{lemma} \label{immediate} Let $U$ be a scheme over $\Z$ and let $k$ be an integer. It follows $\chi(U,k)$ is an integer if and only if Conjecture \ref{sou1} holds for $i=k$. 
\end{lemma}
\begin{proof} The proof is immediate.
\end{proof}

Hence by Lemma \ref{immediate} it follows Conjecture \ref{sou1} is equivalent to the following conjecture:

\begin{conjecture} \label{eulerconj} Let $X$ be a quasi projective scheme of finite type over $\Z$ and let $k$ be an integer. The Euler characteristic $\chi(X,k)$ is an integer.
\end{conjecture}

\begin{example} Conjecture \ref{souleconjecture} for Dedekind L-functions.\end{example}
If $S:=\Spec(\O_K)$ with $K$ an algebraic number field, it follows \ref{sou1}, \ref{sou2} and \ref{sou3} holds by the work of Borel \cite{Borel}.

\begin{example} Conjecture \ref{souleconjecture} for finite fields.\end{example}

Let $k$ be a finite field. It follows $\K'_m(k)_{\Q}=0$ hence $\K'_m(k)_{(j)}=0$ for all integers $m,j$, and it follows \ref{sou1} holds for $S:=\Spec(k)$. One also checks
\ref{sou3} holds for $S$.

Note: In the case when $X$ is a regular scheme of dimension $D$ it follows there is an equality of groups 

\[ \K’_m(X)_{(i)} \cong \K_m(X)_{\Q}^{(D-i)}\]

 $\K_m(X)$ is the K-theory of the category $P(X)$ of locally trivial finite rank $\O_X$-modules.

\begin{lemma} \label{lfunction} Let $X$ be of finite type over $S$ with $X=U \cup V$ a disjoint union of two subschemes $U,V$. It follows $\L(X,s)=\L(U,s)\L(V,s)$. If $U\subseteq X$ is an open subscheme
with $Z:=X-U$ it follows $\L(X,s)=\L(U,s)\L(Z,s)$. Assume $X,Y$ are schemes of finite  type over $S$ such that for any closed point $t\in S$ there is an isomorphism $X_t\cong Y_t$ of fibers.
It follows there is an equality of L-functions $\L(X,s)=\L(Y,s)$. There is an equality $\L(\A^d_X,s)=\L(X,s-d)$. More generally if $E$ is a vector bundle of rank $d$ on X it follows $\L(E,s)=\L(X,s-d)$.
\end{lemma}
\begin{proof} Let $Z \subseteq X$ be a closed subscheme and let $x \in Z$ be a closed point. We let $N_Z(x):= \# \kappa_Z(x)$ denote the number of elements in the residue field $\kappa_Z(x)$.  Here $\kappa_Z(x)$ indicates we view $x$ as a closed point in $Z$.
It follows there is an equality $N_Z(x)=N_X(x)$.

Assume we may write $X$ as a disjoint union $X = U \cup V$. It follows $X^{cl} = U^{cl} \cup V^{cl}$. We get

\[ \L(X,s)= \prod_{x \in X^{cl}}\frac{1}{1-N_X(x)^{-s}} = \prod_{x \in U^{cl}}\frac{1}{1-N_U(x)^{-s}} \prod_{x \in V^{cl}}\frac{1}{1-N_V(x)^{-s}}=\]

\[ \L(U,s)\L(V,s)\]
and the first claim follows. We moreover get

\[ \L(X,s)=\prod_{x \in X^{cl}}\frac{1}{1-N_X(x)^{-s}}= \prod_{t\in S^{cl}} \prod_{x \in X_t^{cl}} \frac{1}{1-N_{X_t}(x)^{-s} }=\]
\[ \prod_{t\in S^{cl}} \L(X_t,s)=\prod_{t\in S^{cl}} \L(Y_t,s)=\L(Y,s) .\]

By Exercise 5.3 and 5.4 in Appendix C in \cite{hartshorne} we get the following: If $T$ is a scheme of finite type over a finite field $k$ with $q$ elements and $Z(T,t)$ is the Weil zeta function
of $T$, then $Z(T\times_k \A^d_k, t)=Z(T,q^dt)$. Moreover $\L(T,s)=Z(T, q^{-s})$. We get the following: If $x\in X$ is a closed point and $T:=\Spec(\kappa(x))$ and $q:=\#\kappa(x)$, 
it follows the fiber of the map
\[ p: \A^d_X \rightarrow X \]
at $T$ is $\A^d_T$. It follows
\[\L(\A^d_X,s) = \prod_{x \in X^{cl} } \L(\A^d_T,s).\]
We get
\[  \L(\A^d_T, s)=Z(\A^d_T, q^{-s})=Z(T, q^dq^{-s})=Z(T, q^{-(s-d) })=\L(T,s-d).\]
It follows
\[ \L(\A^d_X,s)=\prod_{x\in X^{cl} }\L(\A^d_T, s)=\prod_{x\in X^{cl} }\L(T,s-d)=\L(X,s-d) .\]
Since $\A^d_X$ and $E$ have the same fibers it follows $\L(E,s)=\L(\A^d_X,s)$,  hence $\L(E,s)=\L(X,s-d).$ The Lemma follows.
\end{proof}

\begin{corollary} \label{productformula}  Let $S:=\Spec(\O_K)$ with $K$ a number field and let $X$ be a scheme of finite type over $S$. Let $p:X \rightarrow S$ be the canonical map and let $X_t$ be the fiber of $p$ at $t$ for any point $t\in S$. 
It follows there is an equality of formal products
\[ \L(X,s)=\prod_{t \in S^{cl}} \L(X_t,s) .\]
Let  $t\in S^{cl}$ be a closed point and let $q_t:=N_S(t)$. Let  $Z(X_t,t)$ be the Weil zeta function of $X_t$. It follows there is an equality of formal products
\[  \L(X,s)= \prod_{t\in S^{cl} } Z(X_t,q_t^{-s}) .\]
\end{corollary}
\begin{proof} The proof follows from the proof of Lemma \ref{lfunction} and Lemma \ref{productL}.
\end{proof}

\begin{corollary} \label{corrlfunction} Let $E,F$ be locally trivial $\O_K$-modules of rank $d+1$ and let $S:=\Spec(\O_K)$. It follows $\L(\P(E^*),s)=\L(\P(F^*),s)$. 
Assume $T$ is a regular scheme of finite type over $\O_K$ and $\mathbb{A}^d_T$ is affine d-space over $T$. It follows conjecture \ref{sou3} holds for $T$
if and only if  holds for $\mathbb{A}^d_T$. 
Moreover
\begin{align}
&\label{projl}  \L(\P(E^*),s)=\L(S,s)\L(S,s-1)\cdots \L(S,s-d).
\end{align}
\end{corollary}
\begin{proof}  Since $\P(E^*)$ and $\P(F^*)$ have the same fibers, it follows from Lemma \ref{lfunction} there is an equality $\L(\P(E^*),s)=\L(\P(F^*),s)$

Let $dim(T)=n$. We get
\[  \chi(\A^d_T,k)=\sum_{m\in \Z}(-1)^{m+1}dim_{\Q}(\K'_m(\A^d_T)_{(k)}) =\]
\[ \sum_{m\in \Z}(-1)^{m+1}dim_{\Q}(\K_m(\A^d_T)_{\Q}^{(d+n-k)}) =\]
\[ \sum_{m\in \Z}(-1)^{m+1}dim_{\Q}(\K_m(T)_{\Q}^{(n-(k-d))})=\]
\[ \sum_{m\in \Z}(-1)^{m+1}dim_{\Q}(\K'_m(T)_{(k-d)})=\chi(T,k-d).\]
Hence
\[ \chi(\A^d_T,k)=\chi(T,k-d).\]

Assume $ord_{s=k}(\L(T,s))=\chi(T,k)$. We get
\[ ord_{s=k}(\L(\mathbb{A}^d_T,s))=ord_{s=k}(\L(T,s-d)).\]
Let $t:=s-d$, we get
\[ ord_{t=k-d}(\L(T,t))=\chi(T,k-d)=\chi(\mathbb{A}^d_T, k)\]
hence the conjecture holds for $\mathbb{A}^d_T$. The converse is proved similarly.  
There is an equality of L-functions $\L(\P(E^*),s)=\L(\P(\O_S^{d+1}),s)$ and there is a stratification
\[  \emptyset =X_{-1}\subseteq X_0 \subseteq \cdots \subseteq X_d:=\P(\O_S^{d+1}) \]
with $X_i-X_{i-1}=\A^i_S$. Since $\L(X_i-X_{i-1},s)=\L(\A^i_S,s)=\L(S,s-i)$ Formula \ref{projl} follows using induction.
The Corollary is proved.
\end{proof}

\begin{lemma} \label{mvsequence1} Let $X$ be a scheme of finite type over $S:=\Spec(\O_K)$ with $K$ a number field and let $U,V \subseteq X$ be two open subschemes with $X=U \cup V$. It follows there is a long exact sequence
\begin{align}
\label{mvsequence}  \cdots \rightarrow \K'_m(X)_{\Q}^{(j)}  \rightarrow  \K'_m(U)_{\Q}^{(j)} \oplus \K'_m(V)_{\Q}^{(j)}   \rightarrow \K'_m(U\cap V)_{\Q}^{(j)} \rightarrow \K'_{m-1}(X)_{\Q}^{(j)} \rightarrow \cdots 
\end{align}
of abelian groups for all integers $j\in \Z$.
\end{lemma}
\begin{proof} From \cite{quillen}, Remark 7.3.5 there is a long exact sequence of abelian groups
\begin{align}
\label{mv1sequence}  \cdots \rightarrow \K'_m(X)  \rightarrow  \K'_m(U)  \oplus \K'_m(V)    \rightarrow \K'_m(U\cap V)  \rightarrow \K'_{m-1}(X)  \rightarrow \cdots. 
\end{align}
The proof of the existence of \ref{mv1sequence} uses the long exact localization sequence. When we tensor the long exact localization sequence with the rational number field and take the Adams eigenspace at $j$
the sequence remains exact. This implies the exactness of the  sequence in \ref{mvsequence}.
\end{proof}

\begin{lemma} \label{mvlemma} Let $X$ be a scheme of finite type over $S$ and let $k$ be an integer. Let $U,V \subseteq X$ be open subschemes with the property that $X=U \cup V$ and $\chi(X,k), \chi(U,k)$ and $\chi(V,k)$ are integers.
It follows  $\chi(U\ \cap V,k)$ is an integer and there is an equality
\begin{align}
\label{euleradd}\chi(X,k)=\chi(U,k)+\chi(V,k)-\chi(U\cap V,k)
\end{align}
of integers.
\end{lemma}
\begin{proof} Let $N_1 \leq N_2$ be integers with the property that for all $n\leq N_1$ and $n\geq N_2$ it follows $\K_n'(X)_{\Q}^{(k)}=\K_n'(U)_{\Q}^{(k)}=\K_n'(V)_{\Q}^{(k)}=0$. It follows from the exact sequence
\ref{mvsequence} that $\K_n(U\cap V)_{\Q}^{(k)}=0$. Hence $\chi(U\cap V,k)$ is an integer and there is an equality
\[ \chi(X,k)=\chi(U,k)+\chi(V,k)-\chi(U\cap V,k).\]
The Lemma follows.
\end{proof}

\begin{lemma} \label {eulercover} Assume $X$ is a scheme of finite type over $\Z$ and assume $U_1 \cup \cdots \cup U_s =X$ is an open cover of $X$ with $s \geq 3$. Let $m\in \Z$ be an integer.
Assume
\begin{align}
\label{cond1}&\text{ $\chi(U_{i_1} \cap U_{i_2} \cap \cdots \cap U_{i_l},m)$ is an integer} 
\end{align}
 for all $l=1,..,s$ and $1\leq i_1 < i_2 < \cdots < i_l  \leq s$.  It follows $\chi(X, m)$ is an integer and the following formula holds:
\[  \chi(X,m)=\sum_{l=1}^s (-1)^{l+1}\sum_{1\leq i_1 < \cdots < i_l \leq s} \chi(U_{i_1}\cap U_{i_2} \cap \cdots \cap U_{i_l}, m) .\]
\end{lemma}
\begin{proof} Let $s=3$ and assume $\chi(U_i,m), \chi(U_i \cap U_j,m)$ and $\chi(U_1 \cap U_2 \cap U_3,m)$ are integers. 
There is for every integer $j$ a long exact sequence
\begin{align}
\label{mvsequence2}&  \cdots \rightarrow \K'_m(X)_{\Q}^{(j)}  \rightarrow  \K'_m(U_1 \cup U_2)_{\Q}^{(j)} \oplus \K'_m(U_3)_{\Q}^{(j)}   \rightarrow 
\end{align}
\[  \K'_m(U_1 \cap U_3 \cup U_2 \cap U_3)_{\Q}^{(j)} \rightarrow   \K'_{m-1}(X)_{\Q}^{(j)} \rightarrow,  \]
and by \ref{euleradd} it follows $\chi(U_1 \cup U_2,m)$ and $\chi(U_1 \cap U_3 \cup \chi U_2 \cap U_3,m)$  are integers. From sequence \ref{mvsequence2} it follows $\chi(X,m)$ is an integer and there is an equality
\[ \chi(X,m)= \sum_{l=1}^3 (-1)^{l+1} \sum_{1 \leq i_1 < \cdots <i_l \leq 3} \chi(U_{i_1}\cap \cdots \cap U_{i_l},m) \]
Hence the claim is true for $s=3$. Assume the claim of the Lemma is true for $s\geq 3$. We want to prove the Lemma for a cover with $s+1$ open sets. We may write $X=U_1 \cup \cdots \cup U_{s+1}$ and we get the exact sequence
\begin{align}
\label{mvsequence3}&  \cdots \rightarrow \K'_m(X)_{\Q}^{(j)}  \rightarrow  \K'_m(U_1 \cup \cdots \cup U_s )_{\Q}^{(j)} \oplus \K'_m(U_{s+1})_{\Q}^{(j)}   \rightarrow 
\end{align}
\[  \K'_m(U_1 \cap U_{s+1} \cup \cdots \cup U_s \cap U_{s+1})_{\Q}^{(j)} \rightarrow   \K'_{m-1}(X)_{\Q}^{(j)} \rightarrow. \]
By induction it follows $\chi(U_1\cup \cdots \cup U_s,m)$ and $\chi(U_1 \cap U_{s+1} \cup \cdots \cup U_s \cap U_{s+1},m)$ are integers. It follows from the exactness of sequence \ref{mvsequence3} that $\chi(X,m)$ is an integer.
Again from the exactness of the sequence \ref{mvsequence3} we get the following formula when we pass to the  Euler characteristic:
\[ \chi(X,m)=\chi(U_1\cup \cdots \cup U_s,m)+\chi(U_{s+1},m) -\chi(U_1 \cap U_{s+1} \cup \cdots \cup U_s \cap U_{s+1},m).\]
By the induction hypothesis we get
\[ \chi(X,m)= \sum_{l=1}^s (-1)^{l+1}\sum_{1\leq j_1< \cdots j_l \leq s} \chi(U_{j_1}\cap \cdots U_{j_l},m) + \chi(U_{s+1},m) \]
\[ - \sum_{l'=1}^s (-1)^{l'+1} \sum_{1 \leq i_1 < \cdots i_{l'} \leq s} \chi(U_{i_1}\cap \cdots U_{i_{l'}} \cap U_{s+1}) .\]
We rearrange the sums to get
\[ \chi(X,m)= \chi(U_1,m) + \cdots + \chi(U_{s+1},m) \]
\[- \sum_{1 \leq i_1 < i_2 \leq s }\chi(U_{i_1}\cap U_{i_2},m)- \sum_{1\leq j_2 \leq s} \chi(U_{j_1} \cap U_{s+1},m) +  \]
\[ \sum_{1\leq i_1<i_2<i_3 \leq s}\chi(U_{i_1}\cap U_{i_2} \cap U_{i_3},m) + \sum_{1\leq j_1 < j_2 \leq s} \chi(U_{j_1} \cap U_{j_2} \cap U_{s+1},m) + \cdots \]
\[ +(-1)^{s+1} \chi(U_1 \cap \cdots \cap U_s,m)+  \]
\[ (-1)^{s+1}\sum_{1 \leq j_1<j_2< \cdots < j_{s-1} \leq s} \chi(U_{j_1}\cap \cdots \cap U_{j_{s-1}} \cap U_{s+1},m) \]
\[ +(-1)^{s+2} \chi(U_1 \cap \cdots \cap U_{s+1},m)=\]
\[ \sum_{l=1}^{s+1}(-1)^{l+1}\sum_{1 \leq j_1 < \cdots < j_l \leq s+1}\chi(U_{j_1} \cap \cdots \cap U_{j_l},m),\]
and the Lemma follows by induction.
\end{proof}


We may express the L-function $\L(X,s)$ in terms of an open cover of $X$.

\begin{lemma} \label{lcover} Let $X$ be a scheme of finite type over $\O_K$ with  $K$ a number field and let $U_1,..,U_s$ be an open cover of $X$. Let $I:=(i_1,..,i_l)$ be a set of integers with $1\leq i_1< \cdots < i_l \leq s$ and let $U_I:=U_{i_1}\cap \cdots \cap U_{i_l}$.
The following holds:
\begin{align}
\label{lone}&\L(X,s)= \prod_{l=1}^s ( \prod_{1\leq i_1< \cdots < i_l \leq s}\L(U_I,s))^{(-1)^{l+1}} \\
\label{ltwo}& ord_{s=k}(\L(X,s))=\sum_{l=1}^s(-1)^{l+1}\sum_{1 \leq i_1< \cdots < i_l \leq s}ord_{s=k}(\L(U_I,s)).
\end{align} 
\end{lemma}
\begin{proof} The proof is similar to the proof of Lemma \ref{eulercover} and is left to the reader. If $X= U \cup V$ it follows $\L(X,s)=\L(U,s)\L(V,s)/\L(U\cap V,s)$ and the formulas in \ref{lone} and \ref{ltwo} may be proved using induction.
\end{proof}

\begin{example} A local criteria for Conjecture \ref{sou1} and \ref{sou3} to hold. \end{example}

We get a criteria for Conjecture \ref{sou1} and \ref{sou3} to hold in terms of an open cover.

\begin{theorem} \label{soulecover}  Let $X$ be a scheme of finite type over $\Z$ and let $U_1,..,U_s$ be an open cover of $X$. Let $I:=(i_1,..,i_l)$ be a set of integers. We say the set $I$ satisfies condition $P$ if $1 \leq i_1< \cdots < i_l \leq s$. 
Let $U_I:=U_{i_1}\cap \cdots \cap U_{i_l}$.
Assume Conjecture \ref{sou1} holds for a fixed integer $k$ for all $U_I$ satisfying condition $P$. It follows Conjecture \ref{sou1} holds for $X$ at $k$.
Assume Conjecture \ref{sou3} holds for $U_I$ at an integer $k$ for all $I$ satisfying condition $P$. It follows Conjecture \ref{sou3} holds for $X$ at $k$.
\end{theorem}
\begin{proof} Since Conjecture \ref{sou1} holds for $U_I$ at $k$ for all $i$ satisfying $P$, it follows from Lemma \ref{eulercover} Conjecture \ref{sou1} holds for $X$ at $k$. Let $k$ be any integer. The following holds by Lemma \ref{eulercover}:
\[ \chi(X,k)= \sum_{l=1}^n(-1)^{l+1}\sum_I \chi(U_I,k) .\]
Assume Conjecture \ref{sou3} holds at $k$ for all $U_I$. It follows $ord_{s=k}(\L(U_I,s))=\chi(U_I,k)$ for all $I$ satisfying $P$.
We get from Lemma \ref{lcover}
\[  ord_{s=k}(\L(X,s))=\sum_{l=1}^n(-1)^{l+1}\sum_I ord_{s=k}(\L(U_I,s)) =\]
\[ \sum_{l=1}^n (-1)^{l+1}\sum_I \chi(U_I,k)=\chi(X,k).\]
The Theorem follows.
\end{proof}

\begin{lemma} \label{vbeuler} Let $X$ be a scheme of finite type over $S:=\Spec(\O_K)$ with $K$ a number field and assume $\chi(X,k)$ is an integer. Let $E$ be a rank $l$ vector bundle on $X$ with the property that there is a finite open cover 
$U_1,..,U_s$ of $X$ that trivialize $E$ with the property that Condition \ref{cond1} holds for $U_i$.
It follows $\chi(E,k)$ is an integer with  $\chi(E,k)=\chi(X,k-l)$. It follows Conjecture \ref{sou3} holds for $E$  ifand only if it holds for $X$.
\end{lemma}
\begin{proof} Let  $\pi: E \rightarrow X $ be the structure map and assume there is an open cover $X=U_1 \cup \cdots \cup U_s$ satisfying condition \ref{cond1} with $E_i:=\pi^{-1}(U_i)\cong \A^l_{U_i}$. 
Since $\pi^{-1}(U)\cap \pi^{-1}(V) =\pi^{-1 }(U\cap V)$ it follows for any set of integers $1\leq i_1< \cdots < i_l \leq s$ there is an isomorphism of schemes
\[  E_{i_1}\cap \cdots \cap E_{i_l}\cong \A^l_{U_{i_1}\cap \cdots \cap U_{i_l} }.\]
It follows 
\[ \chi( E_{i_1}\cap \cdots \cap E_{i_l},k) =\chi( \A^l_{U_{i_1}\cap \cdots \cap U_{i_l}},k)  \]
is an integer for any integer $k$ since $ \chi( \A^l_{U_{i_1}\cap \cdots \cap U_{i_l}} ,k) $ is an integer. 

From Lemma \ref{eulercover} we get the following:
\[ \chi(E,k)=\chi(\cup E_i,k)=\].
\[ \sum_{l=1}^s (-1)^{l+1}\sum_{1 \leq i_1 < \cdots < i_l \leq s} \chi(E_{i_1} \cap \cdots \cap E_{i_l},k)=\]
\[ \sum_{l=1}^s (-1)^{l+1}\sum_{1 \leq i_1 < \cdots < i_l \leq s} \chi\A^l_{U_{i_1} \cap \cdots \cap U_{i_l}},k).\]
Let $\A^l_X$ be the trivial rank $l$ vector bundle with structure map $\gamma:\A^l_X \rightarrow X$ and $F_i:=\gamma^{-1}(U_i)$. We get similarly
\[  \chi(\A^l_X,k)=\chi(\cup F_i,k)=\] \[ \sum_{l=1}^s (-1)^{l+1}\sum_{1 \leq i_1 < \cdots < i_l \leq s} \chi(F_{i_1} \cap \cdots \cap F_{i_l},k)=\]
\[ \sum_{l=1}^s (-1)^{l+1}\sum_{1 \leq i_1 < \cdots < i_l \leq s} \chi\A^l_{U_{i_1} \cap \cdots \cap U_{i_l}},k)=\]
\[ \chi(E,k).\]
The last equality follows since $\A^l_X$ is a trivial bundle and hence the open cover $U_i$ gives a trivialization of $\A^l_X$ as well. Hence $\chi(E,k)=\chi(\A^l_X,k)=\chi(X, k-l)$.
Since $\L(E,s)=\L(\A^l_X,s)=\L(X,s-l)$  the last statement follows from Lemma \ref{corrlfunction}. The Lemma is proved.
\end{proof}

\begin{example} Conjecture \ref{sou1} and \ref{sou3} for vector bundles on $\Spec(\O_K)$. \end{example}

We can prove the Beilinson-Soule vanishing conjecture and the Soule conjecture on $\L(E,s)$ for any finite rank vector bundle $E$ on $\Spec(\O_K)$ using Lemma \ref{vbeuler}. Note that for a number field $K$
it follows the Picard group $\Pic(\O_K)$ is nontrivial in general, hence there is in general an infinite number of non-trivial finite rank projective $\O_K$-modules.

\begin{theorem} \label{vbtheorem} Let $S:=\Spec(\O_K)$ with $K$ a number field and let $U\subseteq S$ be an open subscheme with complement $Z$. It follows Conjecture \ref{sou1} and \ref{sou3} holds for $U$.
Let $E$ be a vector bundle on $S$  of rank $d$. It follows Conjecture \ref{sou1} and \ref{sou3} holds for $E$ and there is an equality $\chi(E,k)=\chi(S,k-d)$.
\end{theorem}
\begin{proof} It is clear Conjecture \ref{sou1} and \ref{sou3} holds for $Z$ since $Z$ is a finite set of closed points. There is a long exact localization sequence
\begin{align}
\label{longex} \cdots \rightarrow \K_{m+1}'(U)_{\Q}^{(j)} \rightarrow \K_m'(Z)_{\Q}^{(j)} \rightarrow \K_m'(S)_{\Q}^{(j)} \rightarrow  \K_m'(U)_{\Q}^{(j)} \rightarrow 
\end{align}
\[ \K_{m-1}(Z)_{\Q}^{(j)} \rightarrow \cdots  \] 
and there are integers $N_1 \leq N_2$ with $\K_n'(Z)_{\Q}^{(j)}=\K_n(S)_{\Q}^{(j)}=0$ for all $n \leq N_1$ or $n \geq N_2$. It follows $\K_n(U)_{\Q}^{(j)}=0$ for all $n \leq N_1$ or $n\geq N_2$, hence $\chi(U,j)$ is an integer.
It follows Conjecture \ref{sou1} holds for $U$. There are equalities $\chi(U,k)=\chi(S,k)-\chi(Z,k)$ and $\L(U,s)=\L(S,s)/\L(Z,s)$ and it follows
\[ ord_{s=k}(\L(U,s))=ord_{s=k}(\L(S,s))-ord_{s=k}(\L(Z,s))= \]
\[ \chi(S,k)-\chi(Z,s)=\chi(U,k), \]
hence conjecture \ref{sou3} holds for $U$. 
Let $\pi: E\rightarrow S$  be the structure map of a rank $d$ vector bundle $E$ with a finite open cover $S=U_1 \cup \cdots \cup U_s$ such that $\pi^{-1}(U_i)\cong  \A^d_{U_i}$. 
Since $U_i \subseteq S$ is an open subscheme
it follows Conjecture \ref{sou1} holds for $U_{i_1} \cap \cdots \cap U_{i_l}$ for all $1\leq i_1<\cdots <i_l \leq s$. It follows from Lemma \ref{vbeuler} that Conjecture \ref{sou1} holds for $E$ and there is an equality $\chi(E,k)=\chi(S,k-d)$.
The Lemma follows.
\end{proof}

\begin{example} The Euler characteristic of a projective bundle on $\Spec(\O_K)$. \end{example}

Let $S:=\Spec(\O_K)$ with $K$ a number field and let $E$ be a finitely generated and projective $\O_K$-module of rank $d+1$. Let $\P(E^*)$ be the projective space bundle of $E$ with structure map
$\pi: \P(E^*) \rightarrow S$ and let $U_1,..,U_s$ be an open cover of $S$ with $\pi^{-1}(U_i):=E_i\cong \P^d_{U_i}$ giving a local trivialization of $\P(E^*)$.

\begin{lemma} \label{intlemma} There is for any set of integers  $1 \leq i_1 < \cdots < i_l \leq s$ an isomorphism of open sub-schemes of $\P(E^*)$
\[  E_{i_1}\cap \cdots \cap E_{i_l}\cong \P^d_{U_{i_1}\cap \cdots \cap U_{i_l} }.\]
\end{lemma}
\begin{proof} The proof is immediate since $U_i$ is a local trivialization of $\P(E^*)$ and since $\pi^{-1}(U)\cap \pi^{-1}(V)=\pi^{-1}(U \cap V)$ where $\pi$ is the structure map of $\P(E^*)$..
\end{proof}

\begin{lemma} \label{pbundint} Let $U$ be a scheme of finite type over $\Z$ and assume $\chi(U,k)$ is an integer for any integer $k$. Let $\P^d_U$ be projective space over $U$. It follows $\chi(\P^d_U,k)$ is an integer.
\end{lemma}
\begin{proof} By definition $\P^d_U:=\Proj(\O_U[x_0,..,x_d])$ is the relative proj of the sheaf of commutative $\O_U$-algebras $\O_U[ x_0,..,x_d]$
Assume $d=1$. We get $D(x_0)\cong \A^1_U$ and $V(x_0)\cong U$ and it follows from the long exact localization sequence that $\chi(\P^1_U,k)$ is an integer. The Lemma follows by induction.
\end{proof}

\begin{proposition} \label{eulerprojective} Let $F:=\O_S^{d+1}$ be the trivial rank $d+1$ bundle on $S$ and let $\P(F^*)\cong \P^d_S$ be associated projective space bundle.
It follows $\chi(\P(E^*),k)$ is an integer for any integer $k$ and  there is an equality of Euler characteristics $\chi(\P(E^*),k)=\chi(\P^d_S,k)$ for any integer $k$.
\end{proposition}
\begin{proof} Let $U_1,..,U_s$ be a local trivialization of $\P(E^*)$ with $E_i:=\pi^{-1}(U_i)$ and $\pi: \P(E^*) \rightarrow S$ the structure map.
It follows from Lemma \ref{intlemma} $\cup E_i=\P(E^*)$ is an open cover with isomorphisms
\[  E_{i_1}\cap \cdots \cap E_{i_l}\cong \P^d_{U_{i_1}\cap \cdots \cap U_{i_l} } .\]
From Lemma \ref{pbundint}  it follows 
\[ \chi(E_{i_1} \cap \cdots \cap E_{i_l},k)=\chi (\P^d_{U_{i_1}\cap \cdots \cap U_{i_l} },k)\]
is an integer since $U_{i_1}\cap \cdots \cap U_{i_l} \subseteq S$ is an open subscheme.
From Lemma \ref{eulercover} we get
\[  \chi(\P(E^*),k) = \chi( \cup E_i,k)=\]
\[  \sum_{l=1}^s (-1)^{l+1}\sum_{1 \leq i_1 < \cdots < i_l \leq s} \chi(E_{i_1} \cap \cdots \cap E_{i_l},k)=\]
\[ \sum_{l=1}^s (-1)^{l+1}\sum_{1 \leq i_1 < \cdots < i_l \leq s} \chi(\P^d_{U_{i_1} \cap \cdots \cap U_{i_l}},k).\]
It follows $\chi(\P(E^*),k)$ is an integer for any integer $k$.

Let $\P^d_S:=\P(\O_S^{d+1})$ be the trivial rank $l$ vector bundle with structure map $\gamma:\P^d_S \rightarrow S$ and $F_i:=\gamma^{-1}(U_i)$. We get similarly
\[  \chi(\P^d_S,k)=\chi(\cup F_i,k)=\] \[ \sum_{l=1}^s (-1)^{l+1}\sum_{1 \leq i_1 < \cdots < i_l \leq s} \chi(F_{i_1} \cap \cdots \cap F_{i_l},k)=\]
\[ \sum_{l=1}^s (-1)^{l+1}\sum_{1 \leq i_1 < \cdots < i_l \leq s} \chi(\P^d_{U_{i_1} \cap \cdots \cap U_{i_l}},k)=\]
\[ \chi(\P(E^*),k).\]
The Proposition follows.
\end{proof}

\begin{corollary}\label{eulerprojequal} If $E,G$ are arbitrary rank $d+1$ projective $\O_S$-modules it follows there is an equality
\[ \chi(\P(E^*),k)=\chi(\P(G^*),k) \]
for any integer $k$.
\end{corollary}
\begin{proof} From Proposition \ref{eulerprojective} it follows $\chi(\P(E^*),k)$ and $\chi(\P(G^*),k)$ are integers for any integer $k$. We moreover get an equality  
\[ \chi(\P(E^*),k)=\chi(\P^d_S,k)=\chi(\P(G^*),k) \]
and the Corollary follows.
\end{proof}

\begin{example} \label{cohomological}  A cohomological description of the local L-factors of $\L(X,s)$.\end{example}

Assue $\pi: X \rightarrow S$ is a scheme of finite type over $S$ where $S:=\Spec(\O_K)$ and $K$ is a number field. Assume $X_y:=\pi^{-1}(y)$ is smooth and projective for any closed point $y\in S$.
We may express the L-function $\L(X_y,s)$ of the fiber $X_y$ in terms of the Weil zeta function $Z(X_y, t)$ as follows:
\[  \L(X_y,s)= Z(X_y, q_y^{-s}) \]
where $\kappa(y)=\mathbb{F}_q$ and $q_y=p^n$ with $p>0$ a prime. If $X$ is a smooth projective scheme of finite type over $\kappa(y)$, it is well known the function $Z(X, t)$ is a rational function 
\[  Z(X, t)=\frac{P_1(t)\cdots P_{2n-1}(t)}{P_0(t)\cdots P_{2n}(t)} \]
where $n:=dim(X)$. There is moreover a determinantal formula
\begin{align}
&\label{det}  P_i(t)=\operatorname{det}(1-f^*t; \H^i(\overline{X}, \mathbb{Q}_l) ) 
\end{align}
where $f:\overline{X}\rightarrow \overline{X}$ is the Frobenius morphism and $\H^i(\overline{X}, \mathbb{Q}_l)$ is l-adic etale cohomology (see \cite{dwork}, \cite{hartshorne}).
From Corollary \ref{productformula} it follows the global L-function $\L(X,s)$ is calculated by the l-adic etale cohomology groups $\H^*(\overline{X_y},\mathbb{Q}_{l_y})$ for varying primes $l_y\neq char(\kappa(y))$ via the formula
\[  \L(X,s)=\prod_{y\in S^{cl} }Z( X_y, q_y^{-s}) \]
and Formula \ref{det}. The determinantal formula \ref{det} may be proved using other p-adic cohomology theories (rigid cohomology, cristalline cohomology, prismatic cohomology etc.).
The rationality of the Weil zeta function $Z(X, t)$ was first proved by Dwork in \cite{dwork} in 1960 using p-adic methods.

\begin{corollary} \label{integer} Let $U$ be a scheme over $\Z$ with $\chi(U,k)$ an integer for all $k\in \Z$. Let $E$ be a locally trivial $\O_U$-module of rank $d+1$.
It follows 
\[ \chi(\P(E^*), j) = \sum_{i=0}^d \chi(U, k-i) .\]
Hence it follows $\chi(\P(E^*),k)$ is an integer for all $k\in \Z$.
\end{corollary}
\begin{proof} Since $\P(E^*)$ and $\P^d_U$ have the same fibers it follows $\chi(\P(E^*),k) =\chi(\P^d_U,k)$. By induction it follows
 \[ \chi(\P^d_U, j) = \sum_{i=0}^d \chi(U, k-i) \]
and the Corollary follows since $\chi(U,k)$ is an integer for all integers $k$.
\end{proof}



The following Lemma is by some authors referred to as the \emph{Jouanolou trick}:

\begin{lemma} \label{jouanoloulemma} Let $T:=\Spec(B)$ be an affine scheme of finite type over $\Z$ and let $X\subseteq \P^n_T$ be a quasi projective scheme over $T$. It follows there is an affine scheme $W:=\Spec(B)$
and a surjective map $\pi: W \rightarrow X$ where the fibers of $\pi$ is affine $l$-space $\A^l$.
\end{lemma}
\begin{proof} This is proved in \cite{jouanolou}, Lemma 1.5. 
\end{proof}

The affine $\A^l$-fibration $W$ constructed in Lemma \ref{jouanoloulemma} is an \emph{affine torsor for $U$}.

Note that if $U\subseteq \P^n_{\Z}$ is a quasi projective scheme and $\pi:W \rightarrow U$ is an affine torsor with fiber $\A^l$ constructed in Lemma \ref{jouanoloulemma} it follows 
$\L(W,s)=\L(\A^l_U,s)$ since $W$ and $\A^l_U$ have the same fibers. By construction there is an isomorphism 
\[ \pi^*:\K'_m(U) \cong \K'_m(W)  \]
of abelian groups inducing an isomorphism  
\begin{align}
&\label{homotopy}  \pi^*_{(j-l)}: \K'_m(U)_{(j-l)} \cong \K'_m(W)_{(j)}  
\end{align}
for all integers $j$.  Since $W$ has fibers $\A^l$ it follows $dim(W)=d+l$ where $d:=dim(U)$.
Hence we get the following result:

\begin{lemma} \label{jouanolouchi} Let $U\subseteq \P^n_{\Z}$ be a quasi projective scheme and let $p:W\rightarrow U$ be the torsor constructed in Lemma \ref{jouanoloulemma} with fiber $\A^l$. 
It follows $\L(W,s)=\L(\A^l_U,s)$ and $\chi(W,j)=\chi(U,j-l)$ for all integers $j$.
\end{lemma}
\begin{proof} Since $W$ and $\A^l_U$ have the same fibers it follows $\L(W,s)=\L(\A^l_U,s)$ is an equality of L-functions. By Formula \ref{homotopy} we get an equality

\[ \chi(U,j-l):=\sum_{m\in \Z} (-1)^{m+1}dim_{\Q}(\K'_m(U)_{(j-l)}) =\]

\[   \sum_{m\in \Z} (-1)^{m+1}dim_{\Q}(\K'_m(W)_{(j)})= \chi(W,j) ,\]
hence $\chi(U,j-l)=\chi(W,j)$ and the Lemma follows.
\end{proof}

\begin{lemma}\label{indlemma}  Let $X$ be a scheme of finite type over $\Z$ and let $Z\subseteq X$ be a closed subscheme with open complement $U:=X-Z$.
If conjecture \ref{sou1} and \ref{sou3} holds  for $Z$ and $U$ it follows conjecture \ref{sou1} and \ref{sou3} holds for $X$. There is for all integers $k\in \Z$ an equality
\[  \chi(X,k)=\chi(U,k)+\chi(Z,k).\]
\end{lemma}
\begin{proof} Assume $\K'_m(Z)_{(j)}=\K'_m(U)_{(j)}=0$ for almost all $m$. There is a long exact localization sequence
\[ \cdots \rightarrow \K'_m(Z)_{(j)} \rightarrow \K'_m(X)_{(j)} \rightarrow \K'_m(U)_{(j)} \rightarrow \] 
\[ \rightarrow \K'_{m-1}(Z)_{(j)} \rightarrow \K'_{m-1}(X)_{(j)} \rightarrow \K'_{m-1}(U)_{(j)} \rightarrow \cdots \]
hence there are integers $m_1 \leq m_2$ with the following properties: For all integers $m$ with $m\leq m_1$ or $m_2 \leq m$ it follows $\K'_m(Z)_{(j)}=\K'_m(U)_{(j)}=0$.
It follows by the long exact localization sequence that  $\K'_m(X)_{(j)}=0$ for all $m \leq m_1$ and $ m_2\leq m$, hence Conjecture \ref{sou1} holds for $X$.
If Conjecture \ref{sou3} holds for $Z$ and $U$ we get the following:  $\L(X,s)=\L(Z,s)\L(U,s)$. We get
\[ ord_{s=k}(\L(X,s))=ord_{s=k}(\L(Z,s))+ord_{s=k}(\L(U,s))= \]
\[ \chi(Z,k)+\chi(U,k)=\chi(X,k) \]
since the Euler characteristic is additive with respect to $Z,U$, hence Conjecture \ref{sou3} holds for $X$. The Lemma follows.
\end{proof}

\begin{example} Conjectures \ref{sou1} and \ref{sou3} for affine fibrations over $\O_K$. \end{example}

We may generalize the result on the Soule conjecture for finite rank vector bundles to coherent sheaves on $\Spec(\O_K)$ using Lemma \ref{indlemma}. If $E$ is a finitely generated $\O_K$-module there is an open 
subscheme $U \subseteq S:=\Spec(\O_K)$ (possibly empty) where $E$ is locally trivial of finite rank  (see \cite{matsumura} Theorem 4.10). The closed complement $Z:=S-U=\{s_1,..,s_l\}$ is a finite set of closed points.
We may consider the affine fibration $\pi: \mathbb{V}(E^*) \rightarrow S$ where $\mathbb{V}(E^*):= \Spec(\operatorname{Sym}_{\O_K}(E^*))$

\begin{theorem} \label{affinefib} Conjecture \ref{sou1} and \ref{sou3} holds for $\mathbb{V}(E^*)$.
\end{theorem}
\begin{proof} Using a method similar to the metod in the proof of  Theorem \ref{vbtheorem} it follows Conjecture \ref{sou1} holds for $V:=\pi^{-1}(U)$ since $V$  is a finite rank vector bundle on $U$. 
Since $Z$ is a finite set of closed points and $\pi^{-1}(Z)$ is a disjoint union of affine spaces over finite fields, it follows Conjecture \ref{sou1} holds for $\pi^{-1}(Z)$. It follows from Lemma \ref{indlemma} that Conjecture \ref{sou1} holds for $\mathbb{V}(E^*)$.
A similar reasoning proves that conjecture \ref{sou3} holds for $\mathbb{V}(E^*)$ and the Theorem follows.
\end{proof}

We may reduce the study of Conjecture \ref{sou1} and \ref{sou3} to the study of affine regular schemes of finite type over $\Z$, with a systematic use of localization, induction on dimension and the Jouanolou trick from Lemma \ref{jouanoloulemma}:

\begin{theorem} \label{mainreduction} Assume Conjecture \ref{sou1} and \ref{sou3} holds for any affine regular scheme of finite type over $\Z$. It follows 
Conjecture \ref{sou1} and \ref{sou3} holds for any quasi projective scheme $U$ of finite type over $\Z$. 
\end{theorem}
\begin{proof} One first proves using induction, the long exact localization sequence and Jouanolous trick that Conjecture \ref{sou1} holds for any affine scheme $S:=\Spec(A)$
of finite type over $\Z$. Then again using Jouanolous trick, one proves Conjecture \ref{sou1} holds for any quasi projective scheme $U\subseteq \P^n_{\Z}$ of finite type over $\Z$. 

Assume Conjecture \ref{sou3} holds for all affine regular schemes $S:=\Spec(A)$ of finite type over $\Z$. Let $dim(S)=1$. It follows the singular subscheme $S_s \subseteq S$ is a finite set of closed points
with finite residue fields and Conjecture \ref{sou3} holds for $S_s$. We use here the fact that the K-theory of a scheme $X$ is the same as the K-theory of the associated reduced scheme $X_{red}$. The singular scheme $S_s$ may be non-reduced but we can pass to the reduced scheme associated to $S_s$.
Let $U:=S-S_s$. It follows $U\subseteq \P^n_{\Z}$ is a quasi projective regular scheme and hence there is a affine torsor
$p:W\rightarrow U$ with fibers affine $l$-space $\mathbb{A}^l$. It follows since $W$ is an $\A^l$-fibration that $W$ has the same fibers as relative affine space $\mathbb{A}^l_U$ over $U$. Hence by Lemma \ref{corrlfunction} it follows
there is an equality of L-functions 
\[ \L(W,s)=\L(\mathbb{A}^l_U,s).\]
Since $W$ is an affine regular scheme of finite type over $\Z$ it follows Conjecture \ref{sou3} holds for $W$. We get by Lemma \ref{jouanolouchi}
\[  ord_{s=k}(\L(\mathbb{A}^l_U,s))=ord_{s=k}(\L(W,s))=\chi(W,k)= \]
\[  \chi(U,k-l)=\chi(\A^l_U,k) .\]
Hence Conjecture \ref{sou3} holds for $\A^l_U$.
By Lemma \ref{corrlfunction} since Conjecture \ref{sou3} holds for $\mathbb{A}^l_U$ it holds for $U$. Hence Conjecture \ref{sou3} holds for $S_s$ and $U$ and hence it holds for $S$. By induction
on the dimension it follows \ref{sou3} holds for any affine scheme $S$ of finite type over $\Z$.

Assume $U \subseteq \P^n_{\Z}$ is a quasi projective scheme and let $p:W\rightarrow U$ be an affine torsor
with $W:=\Spec(B)$ where $B$ is a finitely generated $\Z$-algebra. It follows by assumption \ref{sou3} holds for $W$.  By the same argument as above it follows \ref{sou3} holds for $\A^l_U$ and again by 
Lemma \ref{corrlfunction} it follows \ref{sou3} holds for $U$. The Theorem follows.
\end{proof}

Note: A result similar to Theorem \ref{mainreduction} for Conjecture \ref{sou2} is mentioned in \Soule's original paper \cite{Soule} in Example 2.4. 
Theorem \ref{mainreduction} is obtained using slightly different techniques in \cite{Kahn}, Lemma 43. The proof of the theorem is not difficult, but I prefer to call it a Theorem, since it is a significant reduction.
The Jouanolou-Thomason trick in its most general form is a generalization of Lemma \ref{jouanoloulemma} to the case of a quasi compact quasi separated scheme with an ample family of line bundles. Conjecture \ref{sou1} and \ref{sou3}  is stated for quasi projective schemes of finite type over  $\Z$.


\begin{example} Conjecture \ref{souleconjecture} for Abelian schemes.\end{example}

Let $A \subseteq \P^n_{T}$ is a projective abelian scheme of finite type over $T:=\Spec(B)$, where $K$ is an algebraic number field and $B$ a finitely generated and regular $\O_K$-algebra.
If Conjecture \ref{sou1} and \ref{sou3} holds for all affine regular schemes $\Spec(A)$ of finite type over $\Z$, it follows from Theorem \ref{mainreduction} Conjecture 
\ref{sou1} and \ref{sou3} holds for any abelian scheme $A\subseteq \P^n_T$. Hence we have reduced the study of the Birch and Swinnerton-Dyer conjecture for abelian schemes to the study of affine regular schemes $\Spec(A)$ of finite type over $\Z$.

\begin{example} Algebraic K-theory for an affine regular scheme of finite type over $\Z$.\end{example}

Let $S:=\Spec(A)$ where $A$ is a finitely generated and regular $\Z$-algebra. It follows from \cite{weibel}, Section IV, 1.16.1 there is an embedding 
\begin{align}
&\label{hurewicz} \K_*(S)\otimes \Q \subseteq \H_*(\operatorname{GL}(A), \Q) 
\end{align}
where $\operatorname{GL}(A)$ is the infinite general linear group of $A$. The embedding in \ref{hurewicz} realize $\K_*(S)\otimes \Q$ as the primitive elements in the Hopf algebra 
$\H_*(\operatorname{GL}(A),\Q)$. There are Adams operators on $\H_*(\operatorname{GL}(A),\Q)$ inducing the classical Adams operators on $\K_*(S)\otimes \Q$, hence the
weight spaces $\K'_m(S)_{(j)}$ may be constructed using the Hopf algebra structure on $\H_*(\operatorname{GL}(A),\Q)$.
In the paper \cite{Borel} Borel calculates the K-groups $\K_*(\O_K)\otimes \Q$ for any algebraic number field $K$ using the embedding \ref{hurewicz}. This is Theorem \ref{boreltheorem}.

If $\mathfrak{gl}(A)$ is  the Lie algebra of infinite matrices with coefficients in $A$, and $A$ is a $\Q$-algebra, it follows by the Loday-Quillen-Tsygan Theorem (see \cite{weibel0}, Theorem 9.10.10) there is an isomorphism
\[ \operatorname{Prim}_n(\H_*(\mathfrak{gl}(A),\Q)) \cong \operatorname{HC}_{n-1}(A), \]
where $\operatorname{HC}_{n-1}(A)$ is cyclic homology of $A$. Hence in this case there is an explicit formula for the space of primitive elements in terms of cyclic homology. One may ask for a "similar"
explicit formula for the space of primitive elements in $\H_*(\operatorname{GL}(A),\Q)$. The space
\[ \operatorname{Prim}(\H_*(\mathfrak{gl}(A),\Q)) \]
is sometimes referred to as the \emph{additive K-theory} of $A$.

\begin{example} Some speculations on a cohomological formulation of the \Soule conjecture.\end{example}

Let $S:=\Spec(A)$ and let $T:=\Spec(\O_K)$. There are operators 
\[ \phi^{m,i}: \H_*(\operatorname{GL}(A),\Q) \rightarrow \H_*(\operatorname{GL}(A),\Q) \]
with the property that the induced morphism
\[ \phi^{m,i}: \operatorname{Prim}_m(\H_*(\operatorname{GL}(A),\Q)) \rightarrow \operatorname{Prim}_m(\H_*(\operatorname{GL}(A),\Q)) \]
has the following property: Let $E(m,i)$ be the set of elements $x$ with $\phi^{m,i}(x)=x$. It follows there is an equality
$E(m,i)=\K'_m(S)_{(i)}$. Hence we may define the Euler characteristic $\chi(S,i)$ using the homology $\H_*(\operatorname{GL}(A),\Q)$ of the infinite general linear group:
\begin{align}
&\label{eulerhom} \chi(S,i):=\sum_{m\in \Z}(-1)^{m+1}dim_{\Q}(E(m,i)).
\end{align}
In \ref{eulerhom} we have not used algebraic K-theory $\K'_m(S)_{(i)}$ to define $\chi(S,i)$. 
By \ref{cohomological} we may define the local L-factors $\L(S_t, s)$ for any closed point $t\in T$ using a p-adic cohomology theory (or l-adic etale cohomology when $S_t$ is smooth and projective).
It follows by the product formula  
\[ \L(S,s)=\prod_{t \in T^{cl}} \L(S_t,s) \]
that the L-function $\L(S,s)$ has a "cohomological description". Hence the \Soule conjecture may be stated as follows: There is for every integer $i$ an equality
\[ \sum_{m\in \Z}(-1)^{m+1}dim_{\Q}(E(m,i))=ord_{s=i}(\L(S,s)) .\]
Hence we may argue that the \Soule conjecture can be formulated "using cohomology and homology" groups associated to the affine scheme $S$. There are precise conjectures
on the existence of an "arithmetic cohomology theory" that simultaneously generalize the algebraic K-theory of $S$ (or homology of the infinite general linear group $\operatorname{GL}(A)$)
and p-adic cohomology of the fibers $S_t$ for all closed points $t \in T$, and what properties such a theory must have in order to prove the \Soule conjecture (see \cite{deninger})

\section{Conjectures on L-functions for schemes with a generalized cellular decomposition.}

In this section we prove the Beilinson-\Soule vanishing conjecture \ref{sou1} and \Soule conjecture \ref{sou3} for any scheme equipped with a generalized cellular decomposition (see Lemma \ref{lemmagen}).
We also prove the conjectures for any projective bundle $\P(E^*)$ and the partial  flag bundle $\mathbb{F}(N,\E)$ of a coherent $\O_S$-module $\E$,  with $S:=\Spec(\O_K)$. Here $K$ is any number field and $\O_K$ is the ring of integers in $K$
(see Corollary \ref{s1flag}, Theorem \ref{s3flag} and Corollary \ref{s13corollary}). Hence for each number field $K$ we get an infinite number of non-trivial examples of schemes in any dimension where Conjecture \ref{sou1} and \ref{sou3} hold (see Example \ref{nontrivial}).


\begin{example} \label{pdbundle} Conjecture \ref{sou1} and \ref{sou3} for the complete flag bundle $\mathbb{F}(E)$. \end{example}

In the following theorem we use Corollary  \ref{eulerprojequal} to give an explicit and elementary proof of Conjecture \ref{sou1} and \ref{sou3} for any projective bundle $\P(E^*)$ on $\Spec(\O_K)$ with $K$ a number field.

\begin{theorem}  \label{lfunctpbundle} Let $K$ be an algebraic number field and let $S:=\Spec(\O_K)$. Let $\P(E^*)$ be a $\P^d$-bundle on $S$. It follows Conjecture \ref{sou1} and \ref{sou3} holds for $\P(E^*)$.
\end{theorem}
\begin{proof} By induction there is the following result:
\[ \chi(\P(E^*), k)=\sum_{i=0}^d \chi(S,k-i) \]
and since $\chi(S,j)$ is an integer for all integers $j$ it follows $\chi(\P(E^*), k)$ is an integer for all integers $k$. By Lemma \ref{immediate} it follows Conjecture \ref{sou1} holds for $\P(E^*)$.
Let 
\[ \P^n_S:=\Proj(\O_K[x_0,..,x_n]) \]
be projective $n$-space over $\O_K$. 
Let $E$ be a rank $n+1$ projective $\O_K$-module and let $\P(E^*)$ be the $\P^n$-bundle of $E$. It follows by Lemma \ref{corrlfunction} that  $\L(\P(E^*),s)=\L(\P^n_S,s)$. 
By Corollary \ref{eulerprojequal} there is an equality of Euler characteristics
\[ \chi(\P(E^*),j)=\chi(\P^n_S,j),  \] 
hence Conjecture \ref{sou3} holds for $\P(E^*)$ if and only if it holds for $\P^n_S$.

Let $n=1$ and let $\P^1_S:=\Proj(\O_K[x_0,x_1])$. Let $S:=V(x_1)\cong \Spec(\O_K):=S$ and let $D(x_1)=\cong \A^1_S=\Spec(\O_K[\frac{x_0}{x_1}])$. It follows

\[ \chi(\P^1_S,k)=\chi(\A^1_S,k)+\chi(S,k) \]
and
\[ \L(\P^1_S,s)=\L(\A^1_S,s)\L(S,s).\]
Hence
\[ ord_{s=k}(\L(\P^1_S,s))=ord_{s=k}(\L(\A^1_S,s))+ord_{s=k}(\L(S,s))= \]
\[  \chi(\A^1_S,k)+\chi(S,k)=\chi(\P^1_S,k), \]

and it follows \ref{sou3} holds for any $\P^1$-bundle on $S$. Assume the conjecture holds for any $\P^{d-1}$-bundle on $S$ and consider $\P^d_S:=\Proj(\O_K[x_0,..,x_d])$. Let $Z:=V(x_d)$ and let $U:=D(x_d)$. It follows $Z\cong \P^{d-1}_{S}$ and $U\cong \A^d_{S}$. Hence Conjecture \ref{sou3}
holds for $Z$ and $U$. By Lemma \ref{corrlfunction} it follows \ref{sou3} holds for $\P^d_S$ and $\P(E^*)$ for any $E$. The Theorem is proved.
\end{proof}

Note: Theorem \ref{lfunctpbundle} is a generalization of Borel's classical result on $\Spec(\O_K)$ to higher dimensional schemes. The picard group $\operatorname{Pic}(\O_K)$ is a finite nontrivial group in general,
and given any set of elements $\mathcal{L}_i$ for $i=0,..,d$ we get a locally trivial $\O_K$-module $E:=\oplus \mathcal{L}_i$ of rank $d+1$ and a $\P^d$-bundle $\P(E^*)$. 

\begin{corollary} \label{globalpd} Let $U$ be a scheme over $\Z$ where Conjecture \ref{sou1} and \ref{sou3} holds and let $E$ be a rank $d+1$ locally trivial $\O_U$-module. It follows 
Conjecture \ref{sou1} and \ref{sou3} holds for the projective bundle $\P(E^*)$ of $E$.
\end{corollary}
\begin{proof} 
The proof is similar to the proof of Theorem \ref{lfunctpbundle} and is left to the reader.  The Corollary is proved.
\end{proof}

\begin{example} Conjectures \ref{sou1} and \ref{sou3} for projective fibrations on $\O_K$. \end{example}

Let $K$ be a number field with ring of integers $A:=\O_K$ and let $S:=\Spec(A)$. Let $E$ be a finitely generated $A$-module and let $\pi: \P(E^*)\rightarrow S$ be the projective space bundle of $E$. It follows there is an open subscheme
$U\subseteq S$ with $\P(E^*_U)\cong \pi^{-1}(U)$ and where $\P(E^*_U)$  is a projective bundle of finite rank over $U$. 
The closed complement $Z:=S-U$ is a finite set of closed points and $\pi^{-1}(Z)\cong \P(E^*_Z)$ is a finite disjoint union of projective spaces over finite fields.

\begin{theorem} Conjecture \ref{sou1} and \ref{sou3} holds for $\P(E^*)$.
\end{theorem}
\begin{proof} By Corollary \ref{globalpd} it follows Conjectures \ref{sou1} and \ref{sou3} holds for $\P(E^*_U)$ since $\P(E^*_U)$ is a finite rank projective bundle over $U$, and since the conjectures hold for $U$. 
Since $\P(E^*_Z)$ is a finite disjoint union of projective spaces over finite fields, it follows Conjecture \ref{sou1} and \ref{sou3} holds
for $\P(E^*_Z)$. It follows from Lemma \ref{indlemma} Conjecture \ref{sou1} and \ref{sou3} holds for $\P(E^*)$.
\end{proof}

\begin{example} \label{sequencepr} Sequences of projective bundles.\end{example}

Let $U$ be a scheme over $\Z$ and construct $X$ as follows: Let $E_1$ be a locally trivial $\O_U$-module of rank $d_1+1$. Let $X_1:=\P(E_1^*)$. Let $E_2$ be a locally trivial $\O_{X_1}$-module
of rank $d_2+1$ and let $X_2:=\P(E_2^*)$. Continue this process to arrive at a scheme $X:=X_e:=\P(E_e^*)$ with a projection morphism
\begin{align}
&\label{projseq}  \pi: X_e \rightarrow U.
\end{align}

\begin{example} \label{fullflag} An explicit construction of the partial flag bundle $\mathbb{F}(N,E)$.\end{example}

Recall the following construction of the partial flag bundle $\Fl(N,E)$ of a locally free sheaf $E$ using grassmannian bundles. Let $U$ be a scheme over $\Z$ and let $E$ be a locally trivial $\O_U$-module of rank $n$. Let $N:=\{n_1,..,n_l\}$
be a set of positive integers with $\sum_i n_i=n$. Let $\G_1:=\G(n_1,E)$ be the grassmannian bundle of rank $n_1$ subbundles of $E$. There is a tautological rank $n_1$ sub-bundle
\begin{align}
&\label{taut1}  \mathcal{S} \subseteq \pi_1^*E  
\end{align}
where $\pi_1:\G_1 \rightarrow U$ is the projection morphism. 
We get an exact sequence of locally trivial $\O_{\G_1}$-modules
\[ 0 \rightarrow \mathcal{S} \rightarrow \pi_1^*E \rightarrow \mathcal{Q}_2 \rightarrow 0.\]
For any morphism of schemes $f:V\rightarrow U$ there is a canonical isomorphism
\[  V\times_U \G(n_1,E)\cong \G(n_1, f^*E) .\]
As a particular case let $s\in U$ be a point with residue field $\kappa(s)$. We get an inclusion map
\[ i: \Spec(\kappa(s)) \rightarrow U .\]
There is by construction a one to one correspondence between maps of schemes over $U$
\[  g: \Spec(\kappa(s)) \rightarrow  \G_1 \]
and inclusions of $\kappa(s)$-vector spaces
\[  \mathcal{S}(g(s)) \subseteq E(s) \]
where $\mathcal{S}(g(s))$ is the fiber of $\mathcal{S}$ at $g(s)\in \G_1$. The $\kappa(s)$-vector space $\mathcal{S}(g(s))$ has by definition dimension $n_1$.
We get a one-to-one correspondence between the $\kappa(s)$-rational points in the fiber $\pi_1^{-1}(s)$, and subspaces  $W  \subseteq E(s)$ of dimension $n_1$. By functoriality there is an isomorphism
of schemes over $\kappa(s)$
\[ \pi_1^{-1}(s)\cong \Spec(\kappa(s))\times_U \G(n_1,E) \cong \G(n_1, i^*E)\cong \G(n_1,E(s)) ,\]
where $\G(n_1,E(s))$ is the classical grassmannian scheme parametrizing $n_1$-dimensional subspaces of the fiber $E(s)$ of $E$ at $s$. Hence $(\G_1,\pi_1)$ is a fibration over $U$ with fibers grassmannian schemes.

Let $\G_2:=\G(n_2,\mathcal{Q}_2)$. We get a canonical projection map $\pi_2:\G_2 \rightarrow U$ with the property that the fiber $\pi_2^{-1}(s)$ is isomorphic to the flag variety $\Fl(n_1,n_2,E(s))$ parametrizing
flags
\[   W_1 \subseteq W_2 \subseteq E(s) \]
of $\kappa(s)$-vector spaces with $dim_{\kappa(s)}(W_i)=n_1+ \cdots +n_i$ for $i=1,2$.
Continue this process to get a scheme $\mathbb{F}(N,E):=\G(n_{l-1}, \mathcal{Q}_{l-1})$ and a projection morphism
\[ \pi: \mathbb{F}(N,E) \rightarrow U.\]
it follows that for any point $s\in U$ it follows the $\kappa(s)$-rational points of the fiber $\pi^{-1}(s)$ corresponds to  flags
\[  0 \neq W_1 \subseteq W_2 \subseteq \cdots \subseteq W_{l-1} \subseteq E(s) \]
with $dim_{\kappa(s)}(W_i)=n_1+\cdots + n_i$ for $i=1,..,l$. It follows there is an isomorphism of schemes
\[ \pi^{-1}(s)\cong \Fl(N,E(s)) ,\]
where $\Fl(N,E(s))$ is the flag scheme of flags of type $N$ in $E(s)$. The scheme $(\mathbb{F}(N,E), \pi)$ is the \emph{flag bundle of $E$ of type $N$}. 
Let $E_n:=\pi^*E$. It follows $E_n$ is a locally trivial $\O$-module on $\mathbb{F}(N,E)$. There is a sequence of locally free sheaves
\begin{align}
&\label{univ} 0 \neq E_1 \subseteq E_2 \subseteq \cdots \subseteq E_{l-1} \subseteq E_n 
\end{align}
on $\mathbb{F}(N,E)$ and $E_i$ is locally trivial of rank $n_i$. The sequence \ref{univ} is the \emph{universal flag} on $\mathbb{F}(N,E)$. 
There is a stratification of closed subschemes
\[ \emptyset = X_{-1} \subseteq X_0 \subseteq \cdots \subseteq X_n:=\Fl(N,E) \]
with the following property: $dim(X_i)=i+dim(U)$, and there is a decompositon $X_i-X_{i+1}=\cup_{j=1,..,n_i} U_{i,j}$ into a finite disjoint union of open subschemes $U_{i,j}\subseteq X_i$ 
with $U_{i,j}\cong \A^i_U$ an isomorphism of schemes over $U$ for $i=1,2,..,n_i$. The construction and basic properties of the partial flag bundle is done in complete generality in \cite{EGA1}.

\begin{example} The complete flag bundle and projective bundles.\end{example}

If $l=n$ and $n_i=1$ for all $i$ it follows $\Fl(N,E)$ is the \emph{complete flag bundle} of $E$. By the above construction we may realize $\Fl(N,E)$ as a "sequence of projective bundles".

\begin{example} The partial flag variety of a vector space over a field. \end{example}

 If $U:=\Spec(k)$ with $k$ a field,  and $E$ an $n$-dimensional $k$-vector space and let $N:=\{n_1,..,n_l\}$ with $\sum_i n_i=n$.  It follows $\mathbb{F}(N,E)$ is the classical 
partial flag variety of $E$ of type $N$, parametrizing flags 
\begin{align}
&\label{flag} 0 \neq W_1 \subseteq W_2 \subseteq \cdots \subseteq W_{l-1} \subseteq E 
\end{align}
in $E$. Here $W_i$ is a $k$-vector subspace of $E$ of dimension $n_1+\cdots + n_i$. This means there is a one-to-one correspondence between the set of $k$-rational points $\Fl(N,E)(k)$ of the flag
variety $\Fl(N,E)$ and the set of flags $\{W_i\}$ of type $N$ in $E$.
If $\SL(E)$ is the special linear group on $E$ and $P\subseteq \SL(E)$ ie the subgroup of elements fixing a flag $\{W_i\}$ in $E$
of type $N$, it follows we may use geometric invariant theory to construct the quotient variety $\SL(E)/P$. It follows $\SL(E)/P$ is canonically isomorphic to the flag variety $\Fl(N,E)$. Hence there is a canonical left action of $\SL(E)$ on $\Fl(N,E)$. 

Hence for the partial flag bundle $\pi: \Fl(N,E)\rightarrow U$ with $E$ a locally trivial $\O_U$-module of rank $n$, it follows the fiber $\pi^{-1}(s)$ may be realized as a quotient
$\SL(E(s))/P(s)$ where $P(s)\subseteq \SL(E(s))$ is a parabolic subgroup.

\begin{corollary} Let $U$ be a scheme over $\Z$ such that  Conjecture  \ref{sou1} and \ref{sou3} holds for $U$ and let $E$ be a locally trivial $\O_U$-module of rank $n$. Let $X_e$ be the scheme constructed in \ref{projseq}. 
It follows Conjecture \ref{sou1} and \ref{sou3} holds for $X_e$. In particular it follows Conjecture \ref{sou1} and \ref{sou3} holds for the full flag bundle $\mathbb{F}(N,E)$ of $E$. 
\end{corollary}
\begin{proof} The first part of the Corollary follows from \ref{globalpd}, \ref{integer} and an induction. The full flag bundle $\mathbb{F}(N,E)$ is by \ref{fullflag} constructed using projective bundles and the Corollary follows.
\end{proof}


\begin{theorem} \label{beisouvanishing} 
Let $S:=\Spec(A)$ where $A$ is a finitely generated and regular $\Z$-algebra. Let $X\subseteq \P^{n}_S$ be a quasi projective regular 
scheme of dimension $d$. If conjecture \ref{sou1} holds for all affine regular schemes of finite type over $\Z$, it follows Conjecture \ref{sou1} holds for $X$.
\end{theorem}
\begin{proof} By Lemma \ref{jouanoloulemma} there is an affine torsor
\[ p:W \rightarrow X \]
with $W:=\Spec(B)$ with $dim(W)=d+l$. The map $p$ induce an isomorphism at K-theory
\[ p_*: \K'_m(X) \rightarrow \K'_m(W) \]
and weight spaces
\[  p_*: \K'_m(X)_{(i)} \rightarrow \K'_m(W)_{(i+l)} .\]
Since $W$ is affine and finite dimensional it follows for a fixed $i+l$ the group $\K'_m(W)_{(i+l)}=0$ for almost all $m$ by assumption.
Hence the same holds for $\K'_m(X)_{(i)}$. The Theorem follows.
\end{proof}

\begin{corollary} Let $A$ be a finitely generated and regular $\Z$-algebra and let $X\subseteq \P^n_S$ be a quasi projective and regular scheme with $S:=\Spec(A)$. 
Assume Conjecture \ref{sou1} holds for all affine regular schemes of finite type over $\Z$. It follows $\chi(X,i)$ is an integer
for all $i\in \Z$.
\end{corollary}
\begin{proof} This follows from Theorem \ref{beisouvanishing}, since in this case $\chi(X,i)$ is a finite sum of integers.
\end{proof}

\begin{example} The projective bundle formula and the Adams operation.\end{example}

In the following we calculate the K-theory of any finite rank  projective bundle on $S:=\Spec(\O_K)$ using the projective bundle formula and Borel's calculation of $\K'_m(\O_K)$.

The \emph{projective bundle formula} says the following. There is a canonical pull back morphism
\[ \pi^*: \K_*(S) \rightarrow \K_*(\P(E^*))  \]
inducing maps
\[ \pi^*:\K_m(S)_{\Q}^{(i)} \rightarrow \K_m(\P(E^*))_{\Q}^{(i)}  \]
and an isomorphism

\begin{align}
&\label{pbundle}\K_*(\P(E^*)) \cong \K_*(S)\otimes_{\K_0(S)} \K_0(\P(E^*)) \cong \K_*(S)\otimes_{\Z} \Z[t]/(t^{d+1}) .
\end{align}
with $t:=1-L$ and $L:=[\O_{\P(E^*)}(-1)]\in \K_0(\P(E^*))$. The Adams operation $\psi^k$ acts as follows:
\[ \psi^k(t):=1-\psi^k(L)=1-L^k.\]
We get for any element $ zt^j\in \K_m(\P(E^*))\cong\K_m(S)\{1,t,..,t^d\}$ the following formula:
\[ \psi^k(zt^j)=\psi^k(z)(1-L^k)^j \in \K_m(\P(E^*)).\]
The isomorphism
\[  \K_m(\P(E^*)) \cong \K_m(S)\{1,t,..,t^d\} \]
is an isomorphism of $\K_0(S)$-modules.
In Theorem \ref{pbundleE} we use formula \ref{pbundle} and Theorem \ref{boreltheorem} to calculate $\K_m(\P(E^*))$ for all integers $m$.

\begin{theorem} \label{pbundleE} Let $\P(E^*)$ be a $\P^d$-bundle on $S$. The following holds:
\begin{align}
&\label{pbund0}\K_0(\P(E^*))_{\Q}\cong \Q^{d+1}\\
&\label{pbund1}\K_m(\P(E^*))_{\Q}\cong 0\text{ for $m=2i, i\neq 0$}\\
&\label{pbund2}\K_m(\P(E^*))_{\Q}\cong \Q^{r_1+r_2}\otimes \Q^{d+1}\text{ for $m\equiv 1\text{ mod }4$}\\
&\label{pbund3}\K_m(\P(E^*)_{\Q}\cong \Q^{r_2}\otimes \Q^{d+1}\text{ for $m \equiv 3\text{ mod }4$}.
\end{align}
\end{theorem}
\begin{proof} The Theorem follows from Theorem \ref{boreltheorem} and the formula \ref{pbundle}.
\end{proof}

\begin{corollary} Let $T$ be a scheme of finite type over $\Z$ with the property that Conjecture \ref{sou2} holds for $T$.
Let $\P(E^*)$ be a $\P^d$-bundle on $T$. It follows Conjecture \ref{sou2} holds for $\P(E^*)$.
In particular it follows Conjecture \ref{sou2} holds for any $\P^d$-bundle on $\O_K$.
\end{corollary}
\begin{proof} By the projective bundle formula there is an isomorphism of abelian groups
\[ \K'_m(\P(E^*)) \cong \K'_m(T)\{1,t,..,t^d\}.\]
Let $R:=\Q[t]/(t^{d+1})$ with $t:=1-L$. It follows $\psi^k$ acts on $R$ as follows: $\psi^k(t)=\psi^k(1-L)=1-L^k$.
Let $v\in \Z$ be an integer and let $R_{(v)}$ denote the vector space of element $x \in R$ with $\psi^k(x)=k^vx$. It follows there is an inclusion of vector spaces over $\Q$:
\[  \K'_m(\P(E^*))_{(j)} \subseteq \oplus_{u+v=j}\K'_m(T)_{(u)}\otimes R_{(v)}  \]
and since by asumption
\[  dim_{\Q}(  \oplus_{u+v=j}\K'_m(T)_{(u)}\otimes R_{(v)}) < \infty \]
for all $m,j$ it follows $dim_{\Q}(\K'_m(\P(E^*))_{(j)}) < \infty $ for all $m,j$ and the Corollary follows.
\end{proof}




The aim of this section is to prove Conjecture \ref{sou1} and \ref{sou3} for all flag bundles $\Fl(N,E)$ on $S:=\Spec(\O_K)$.
Let $T$ be a fixed regular and quasi projective scheme of finite type over $S:=\Spec(\O_K)$ with $K$ a number field and let $X$ be a scheme of finite type over $T$. Assume 
there is a stratification
\begin{align}
&\label{cellular} \emptyset = X_{-1} \subseteq X_0 \subseteq X_1 \subseteq \cdots \subseteq X_n =X 
\end{align}
of $X$ by closed subschemes $X_i \subseteq X$ with $dim(X_i)=i+dim(T)$.

\begin{definition} \label{celdef} We say the stratification $\{X_i\}_{i=0,..,n}$ is a \emph{cellular decomposition} of $X$ if the following holds: For each $i$ there is an isomorphism (as subschemes of $X_i$)
\[ X_i-X_{i-1} =\cup_{j} U_{i,j} \]
where $\cup_j U_{i,j}$ is a finite disjoint union of open subschemes $U_{i,j}\subseteq X_i$, with isomorphisms $f_{i,j}:U_{i,j} \cong \A^i_T$ where $\A^i_T$ is affine $i$-space over $T$. The map $f_{i,j}$
is an isomorphism of schemes over $T$.
\end{definition}

\begin{theorem} \label{specialmain} Let $T$ be a regular quasi projective scheme of finite type over $\O_K$ such that Conjecture \ref{sou1} and \ref{sou3} holds for $T$. Let $X$ be a  
scheme of finite type over $T$ with a cellular decomposition
\[ \emptyset=X_{-1} \subseteq X_0 \subseteq \cdots \subseteq X_{n-1} \subseteq X_n:=X \]
with $X_i-X_{i+1}=\cup_{j=1,..,n_i} \A^i_T$. It follows Conjecture \ref{sou1} holds for $X$. Moreover
\[ ord_{s=k}(\L(X,s))=\chi(X,k),\]
hence Conjecture \ref{sou3} holds for $X$.
\end{theorem}
\begin{proof} The proof is by induction. We will repeatedly use the following Lemma: Let $X$ be a scheme of finite type over $T$ and let $U\subseteq X$ be an open subscheme with $Z:=X-U$. If Conjecture 
\ref{sou1} and \ref{sou3} holds for $U$ and $Z$ it follows \ref{sou1} and \ref{sou3} holds for $X$.

Since $X_0=T$ it follows Conjecture \ref{sou1} holds for $X_0$. Let $X_1-X_0 =\cup_{j} U_{1,j}$ be a finite disjoint union of $X_1-X_0$ into 
affine open subschemes $U_{1,j}\cong \A^1_T$. We get the following calculation:
\[ \chi(X_1-X_0,k)=\chi(\cup_j U_{1,j},k) =\sum_j \chi(U_{1,j},k)=\sum_j \chi(\A^1_T,k) ,\]
and since \ref{sou1} holds for $T$ it holds for $\A^1_T$. Hence $\chi(A^1_T,k)$ is an integer for all integers $k$. It follows the finite sum
\[ \chi(X_1-X_0,k)=\sum_j \chi(U_{1,j},k) \]
is an integer for all integers $k$. Hence Lemma \ref{immediate} implies that Conjecture \ref{sou1} holds for $X_1-X_0$. Since the conjecture holds for $X_0=T$ by assumption, 
if follows \ref{sou1} holds for $X_1$. By induction it follows \ref{sou1} holds for $X=X_n$.

Assume \ref{sou3} holds for $X_0=T$ and let $X_1-X_0=\cup_j U_{1,j}$ a finite disjoint union into affine bundles $U_{1,j}\cong \A^1_T$. We get
\[ ord_{s=k}(\L(X_1-X_0,s))=ord_{s=k}(\prod_j \L(\A^1_T,s))=\sum_j ord_{s=k}(\L(\A^1_T,s))=\]
\[ \sum_j \chi(\A^1_T,k)= \chi(\cup_j U_{1,j},k)=\chi(X_1-X_0,k)\]
hence conjecture \ref{sou3} holds for $X_1-X_0$. It follows Conjecture \ref{sou3} holds for $X_1$. By induction it follows Conjecture \ref{sou3} holds for $X_n=X$ and the Theorem follows.
\end{proof}

\begin{example}  Conjecture \ref{sou1} and \ref{sou3} for flag bundles of coherent $\O_S$-modules. \end{example}

If $A$ is any commutative unital ring and $F:=A^{d+1}$ the free $A$-module of rank $d+1$ it follows $\P(F^*)\cong \P^d_S$ where $S:=\Spec(A)$. The scheme $\P^d_S:=\Proj(A[x_0,..,x_d])$ has a cellular decomposition 
defined in terms of a basis for the free $A$-module $F$. Choosing a basis for $F$ gives rise to a set of generators for the ring $A[x_0,..,x_d]$ that is transcendental over $A$ and this generating set gives rise to a cellular decomposition
of $\P^d_S$. If $E$ is a rank $d+1$ projective $A$-module that is not free, it is not clear how to define a global cellular decomposition for $\P(E^*)$. Given an open set $U$ where $E_U$ trivialize, we get a basis for $E_U$ 
as free $\O_U$-module, and such a basis gives rise to a cellular decomposition of $\P^d_U$. Hence $S$ has an open cover $U_i$ such that $\pi^{-1}(U_i)\subseteq \P(E^*)$ has a cellular decomposition. 

More generally let  $S:=\Spec(\O_K)$ with $K$ a number field and $\E$ a finite rank locally trivial $\O_S$-module with local trivialization $U_1,..,U_s$ and where $\E_{U_i}\cong \O_{U_i}^d$. Let $\Fl(D,\E)$ be the partial flag bundle of $\E$
of type $D$ and let $\pi:\Fl(D,\E) \rightarrow S$ be the projection morphism with $V_i:=\pi^{-1}(U_i)$. Let $I:=(i_1,..,i_l)$ be a set of integers satisfying property $P: 1\leq i_1< \cdots < i_l \leq s$ . 
Let $U_I:=U_{i_1} \cap \cdots \cap U_{i_l}$ and $V_I:=V_{i_1}\cap \cdots \cap V_{i_l}$. By functoriality it follows $V_i=\pi^{-1}(U_i)\cong \Fl(D,\E_{U_i})\cong \Fl(D,\O_{U_i}^d)$. Assume $V_i$ has a cellular decomposition
\[     X(i)_{m+1}\subseteq X(i)_m \subseteq V_i \]
with $X_i-X_{i+1}\cong \cup \A^{n_{ij}}_{U_i}$ a disjoint union of trivial vector bundles over $U_i$. It follows Conjecture \ref{sou1} and \ref{sou3} holds for $V_i$ since it holds for $U_i$ and any trivial finite rank vector bundle $\A^l_{U_i}$.
Moreover for any intersection $V_I=V_{i_1}\cap \cdots \cap V_{i_l}\cong \Fl(D, \E_{U_I})$, where $U_I:=U_{i_1}\cap \cdots \cap U_{i_l}$, it follows Conjecture \ref{sou1} holds since by assumption $\Fl(D,\E_{U_I})$ has a cellular decomposition relative to $U_I$.

\begin{lemma} \label{s1flag} Let $k$ be an integer.  It follows Conjecture \ref{sou1} holds for $\Fl(D,\E)$ at $k$ and the following formula holds:
\begin{align}
\label{flageuler}&\chi(\Fl(D,\E),k)= \sum_{l=1}^s(-1)^{l+1}\sum_{1\leq i_1<\cdots < i_l \leq s} \chi(V_I,k)
\end{align}
\end{lemma}
\begin{proof}By assumption Conjecture \ref{sou1} holds for any intersection $V_I:=V_{i_1}\cap \cdots \cap V_{i_l}$ for $1\leq i_1< \cdots < i_l \leq s$. It follows $\chi(V_I,k)$ is an integer for any integer $k$. The open sets $V_i$ cover $\Fl(D,\E)$ and it follows
from Lemma \ref{eulercover} that Conjecture \ref{sou1} holds for $\Fl(D,\E)$. Again by Lemma \ref{eulercover} there is the formula
\[ \chi(\Fl(D,\E),k)= \sum_{l=1}^s(-1)^{l+1}\sum_{1\leq i_1<\cdots < i_l \leq s} \chi(V_I,k)  .\]
The Lemma is proved.
\end{proof}

\begin{theorem} \label{s3flag} Let $S:=\Spec(\O_K)$ with $K$ a number field and let $\E$ be a finite rank locally trivial $\O_S$-module and let $k$ be an integer. 
It follows Conjecture \ref{sou3} holds for the partial flag bundle $\Fl(D,\E)$ at $k$.
\end{theorem}
\begin{proof}  We use the criteria given in Theorem \ref{soulecover} and a local trivialization of the $\O_S$-module $\E$.
Let $U_1,..,U_s$ be a local trivialization of $\E$. Let $I:=(i_1,..,i_l)$ be a set of integers with $1 \leq i_1< \cdots < i_l \leq s$ and let $U_I:=U_{i_1}\cap \cdots \cap U_{i_l}$.
Let $\pi:\Fl(D,\E) \rightarrow S$ be the projection map and let $V_i:=\pi^{-1}(U_i)$ and let $V_I:=V_{i_1}\cap \cdots \cap V_{i_l}$. It follows from Lemma \ref{lcover}
\[ \L(\Fl(D,\E)= \prod_{l=1}^s(\prod_{1\leq i_1< \cdots < i_l \leq s}\L(V_I,s))^{(-1)^{l+1}} .\]
We get by Lemma \ref{lcover} the following:
\[  ord_{s=k}(\L(\Fl(D,\E))=\]
\[  \sum_{l=1}^s(-1)^{l+1}\sum_{1\leq i_1< \cdots < i_l \leq s}ord_{s=k}(\L(V_I,s)) =\]
\[  \sum_{l=1}^s(-1)^{l+1}\sum_{1\leq i_1< \cdots < i_l \leq s} \chi( V_I,k)) = \chi(\Fl(D,\E),k)\]
since by assumption $ord_{s=k}(\L(V_I,s))=\chi(V_I,k)$. The Theorem follows.
\end{proof}

\begin{corollary} \label{s13corollary} Let $K$ be a number field and let $\O_K$ be it's ring of integers. Let $S:=\Spec(\O_K)$  and let $k$ be an integer. 
Let $\E$ be a coherent $\O_S$-module with an nonempty open subscheme $U\subseteq S$ such that $\E_U$ is a finite rank locally trivial $\O_S$-module. It follows  Conjecture \ref{sou1} and \ref{sou3} holds for $\Fl(D,\E)$ at $k$.
\end{corollary}
\begin{proof} Let $Z:=S-U$ be the closed complement of $U$, which is a finite set of closed points and let $\pi:\Fl(D,\E) \rightarrow S$ be the projection morphism. 
By functoriality it follows $\pi^{-1}(Z):=\Fl(D,\E_Z)$ is a finite disjoint union of partial flag schemes $\Fl(D,\E(s))$ where $s\in S$ is a closed point, $\E(s)$ is a $\kappa(s)$ vector space and $\kappa(s)$ is a finite field.
It follows Conjecture \ref{sou1} and \ref{sou3} holds for $\Fl(D,\E_Z)$ at $k$. Since $\E_U$ is a finite rank locally trivial $\O_S$-module it follows with an argument similar to the one in Corollary \ref{s1flag} and Theorem \ref{s3flag} 
that Conjecture \ref{sou1} and \ref{sou3} holds for $\Fl(D,\E_U):=\pi^{-1}(U)$ at $k$.  Since $\Fl(D,\E)=\pi^{-1}(U)\cup \pi^{-1}(Z)$ it follows from Lemma \ref{indlemma} that Conjecture \ref{sou1} and \ref{sou3} holds for $\Fl(D,\E)$ at $k$. 
The Corollary is proved.
\end{proof}

\begin{example} \label{nontrivial} Non trivial examples for $\O_K$. \end{example}

Let $K$ be an algebraic number field with $S:=\Spec(\O_K)$ and $\operatorname{Pic}(\O_K)$ non-trivial. Let $L_1,..,L_n \in \operatorname{Pic}(\O_K)$ and let $E:=\oplus_i L_i$. It follows $E$ is a non trivial locally trivial $\O_K$-module of rank $n$.
Hence the flag bundle $\pi: \Fl(N,E) \rightarrow S$ is a non-trivial partial flag bundle on $S$, with the property that the Beilinson-\Soule vanishing conjecture and the \Soule conjecture on L-functions holds
for $\Fl(N,E)$.

\begin{example} Generalized cellular decompositions.\end{example}

In this example we prove Theorem \ref{specialmain} for a larger class of schemes: Schemes equipped with a cellular decomposition of type $\{T_i\}$.

\begin{definition} \label{gencel}
Let $T_0, \ldots,T_n$ be schemes of finite type over $S:=\Spec(\O_K)$ and Let $X$ be a scheme  of finite type over $T_n$. Assume there is a stratification
\[  \emptyset = X_{-1} \subseteq X_0 \subseteq X_1 \subseteq \cdots \subseteq X_n =X \]
of $X$, where $X_i\subseteq X$ is a closed subscheme for every $i$ with the following property: For any $i$ it follows $E_i:=X_i-X_{i-1}$ is a vector bundle over $T_i$ with fiber $\A^{d_i}$.
We say $\{X_i\}_{i=0,..,n}$ is a \emph{cellular decomposition of $X$ of type $\{T_i\}$}. We also say \emph{$X$ has a generalized cellular decomposition}.
\end{definition}

Note: It is clear a cellular decomposition is a generalized cellular decomposition: From Definition \ref{celdef} it follows the scheme $X_i-X_{i-1}$ is an affine vector bundle over $T_i$ with fiber $\A^i$, since 
$X_i-X_{i-1}$ has an open cover $X_i-X_{i-1}=\cup_j U_{i,j}$with $U_{i,j}\cong \A^i_T$. Let $T_i:=\cup_j T$ for all $i$. Hence if $\cup_j U_{i,j}=U_{i,1}\cup \cdots \cup U_{i,l}$ it follows
$T_i:= T \cup \cdots \cup T$: The disjoint union of $T$ taken $l$ times. It follows $X_i-X_{i-1}$ is an affine finite rank vector bundle over $T_i$.

\begin{lemma}  \label{lemmagen} Assume $X$ has a generalized cellular decomposition $X_i \subseteq X$ of type $\{T_i\}$ and assume Conjecture \ref{sou1} and \ref{sou3} holds for $T_i$. 
Assume furthermore that each vector bundle $E_i \rightarrow X_i-X_{i-1}$ in the stratification has a finite local trivialization satisfying condition \ref{cond1}.
It follows Conjecture \ref{sou1} and \ref{sou3}
holds for $X$.
\end{lemma}
\begin{proof} Since $X_0:=X_0-X_{-1}=E_0$ is a finite rank affine vector bundle over $T_0$ and Conjecture \ref{sou1} holds for $T_0$ it follows from Lemma \ref{vbeuler} Conjeture \ref{sou1} holds for $E_0:=X_0$.
By definition $X_1-X_0:=E_1$ is a finite rank affine vector bundle over $T_1$. Conjecture \ref{sou1} holds for $T_1$ hence from Lemma \ref{vbeuler} it holds for $E_1$. It follows Conjecture \ref{sou1} holds for $X_1$.
By induction and using Lemma \ref{vbeuler} it follows Conjecture \ref{sou1} holds for $X$. Conjecture \ref{sou3} is proved similarly and the Lemma follows.

Assume Conjecture \ref{sou3} holds for $T$ and let $E$ be a finite rank vector bundle on $T$ of rank $d$. Since $\chi(E,k)=\chi(\A^d_T,k)$ and $\L(E,s)=\L(\A^d_T,s)$ it follows from lemma \ref{vbeuler} Conjecture \ref{sou3}
holds for $T$ if and only if it holds for $E$. Since $E_1:=X_1-X_0$ is a finite rank vector bundle over $T_1$ it follows Conjecture \ref{sou3} holds for $E_1$. Since \ref{sou3} holds for $X_0$ and $X_1-X_0$
it holds for $X_1$. By induction it follows Conjecture \ref{sou3} holds for $X_n=X$ and the Lemma follows. 
\end{proof} 

Let $\pi: X\rightarrow T$ be a scheme of finite type over $T$ with the following property: There is a zero dimensional closed subscheme $S\subseteq X$ with $U:=X-S$ a vector bundle over $T$ of rank $l$.
Since Conjecture \ref{sou1} and \ref{sou3} hold for $T$ it follows by Lemma \ref{lemmagen} Conjecture \ref{sou1} and \ref{sou3} hold for $X$. The scheme $X$ does not neccessarily have a cellular decomposition but it has by definition a generalized cellular decomposition. 

\begin{example} A generalized cellular decomposition for abelian schemes \end{example}

Let $S:=\Spec(A)$ where $A$ is a finitely generated and regular over $\Z$ and let $A\subseteq \P^n_S$ be an abelian scheme over $S$. Let $i:Z \rightarrow A$ be a closed sub-scheme with open complement $j: U\rightarrow A$, 
and consider the localization sequence
\[  \CH^*(Z) \rightarrow \CH^*(A) \rightarrow \CH^*(U) \rightarrow  0 \]
where   $\CH^*(A)$ is the Chow- group of $A$. The Chow-group $\CH^*(A)$ is non-trivial in general and assume $Im(i_*)\neq (0)$ and $\CH^*(U) \neq (0)$. Since $\CH^*(A)$ is highly non-trivial, it follows many closed subshchemes $Z$ have this property.
One want to construct a scheme $T$ of finite type over $\Z$ with the property that there is morphism $\pi: U \rightarrow T$ and such that $U$ is a finite rank vector bundle over $T$. It follows $\L(U,s)=\L(T,s-d)$ and
$\chi(U,k)=\chi(T,k-i)$. Hence the study of the Soule conjecture for $U$ is reduced to the study of the same conjecture for $T$. It is a natural question to ask if there is a generalized cellular decomposition of the abelian scheme $A$.
This is a non-trivial open problem.

For a curve $C$  of genus $g$ over an algebraically closed field $k$  it follows the symmertric product $C(d)$ may be realized as the projective space fibration $\P(\E^*)$ of a coherent $\mathcal{O}_{J(C)}$-module $\E$, where $J(C)$ is the jacobian
of $C $ (see \cite{mattuck} and \cite{schwarzenberger}). It may be the methods introduced in this paper can be used in this study. One has to develop a similar formalism for the Neron model $E*$ of a curve $E$ of genus $g$ over a number field $K$.

\begin{example} Explicit formulas of L-functions and Euler characteristics.\end{example}

We get explicit formulas for the L-function and Euler characteristic for a scheme $X$ with a cellular decomposition of type $\{T_i\}$.

\begin{lemma} \label{cellularT}Let $X$ be a scheme of finite type over $\O_K$ with a cellular decomposition of type $\{T_i\}_{i=0,\ldots, n}$, where $T_i$ satisfy Conjecture \ref{sou1} and \ref{sou3}.
Let $E_i:=X_i-X_{i_1}$ be a rank $d_i$ trivial vector bundle on $T_i$ for $i=0,\ldots ,n$.
It follows
\begin{align}
&\label{LX} \L(X,s)=\prod_{i=0}^n \L(T_i, s-d_i)\\
&\label{chX} \chi(X, k)= \sum_{i=0}^n \chi(T_i, k-d_i).
\end{align}
\end{lemma}
\begin{proof} Let $X$ be a scheme of finite type over $\Z$ with $U\subseteq X$ an open subscheme with complement $Z:=X-U$. Let $E \rightarrow X$ be a vector bundle of rank $d$. Using methods from Lemma \ref{lfunction} and \ref{corrlfunction},  it follows $\L(X,s)=\L(U,s)\L(Z,s)$ and $\L(E,s)=\L(X,s-d)$. Moreover $\chi(X,k)=\chi(U,k)+\chi(Z,k)$ and $\chi(E,k)=\chi(X,k-d)$. Using this, the Lemma follows by induction.
\end{proof}

Let $E$ be a rank $n$ locally trivial $\O_S$-module with $S:=\Spec(\O_K)$, and let $\Fl(N,E)$ be the partial flag bundle of $E$ of type $N$. There is a cellular decomposition 
\begin{align}
&\label{strat} \emptyset =X_{-1} \subseteq X_0 \subseteq \cdots \subseteq X_n:= \Fl(N,E) 
\end{align}
with $X_i-X_{i-1}:=E_i$ a rank $i$ trivial vector bundle on $T_i$ for $i=0,\ldots ,n$. Here $T_i:=\prod_{j=1}^{l_i}S$. Hence $dim(T_i)=l_i$.
We get explicit formulas for the L-function and Euler characteristic of $\Fl(N,E)$:

\begin{lemma} \label{positive} The following holds:
\begin{align}
&\label{LF} \L(\Fl(N,E),s)=\prod_{i=0}^n \L(S, s-d_i)^{l_i} \\
&\label{chF} \chi(X, k)= \sum_{i=0}^n l_i \chi(S, k-d_i).
\end{align}
\end{lemma}
\begin{proof} Since any partial flag bundle $\Fl(N,E)$ has a cellular decomposition of type $\{T_i\}$, the Lemma follows from Lemma \ref{cellularT}.
\end{proof}

\begin{example} An alternative approach using induction.\end{example}

Given a locally trivial finite rank $\O_K$-module $E$ and a flag bundle $\Fl(E)$, we may ask if it is possible to give a proof of Conjecture \ref{sou3} using an induction similar to Example \ref{pdbundle}.
One wants a stratification of closed subschemes
\[ \emptyset = X_{n+1} \subseteq X_n \subseteq \cdots \subseteq X_2 \subseteq X_1 =\Fl(E) \]
with $X_i-X_{i+1}=\cup_{i,j} \A^{d_i}$ is a disjoint union of affine spaces, and where the sub-schemes $X_i$ are flag schemes of dimension smaller than $\Fl(E)$ with the property that Conjecture \ref{sou1} hold for $X_i$. This is done in Example \ref{pdbundle}
for $\P^d$-bundles on $\O_K$. In Example \ref{pbundle} the schemes $X_i$ are projective spaces over $S$  of dimension less than $d$. 

Let $k$ be a field, $E$ an $n$-dimensional vector space over $k$ and let $N:=\{n_1,..,n_l\}$ with $\sum_i n_i =n$. Let $E$ have a flag of $k$-vector spaces
\[ E_1 \subseteq E_2 \subseteq \cdots \subseteq E_l \subseteq E \]
with $dim_k(E_i)=n_1+\cdots + n_i$. Let $\Fl(E)$ be the complete flag variety of $E$. It follows there is a Borel subgroup $B\subseteq \SL(E)$ and an isomorphism $\SL(E)/B\cong \Fl(E)$. There is moreover a parabolic subgroup
$P\subseteq \SL(E)$ with $\SL(E)/P\cong \Fl(N,E)$, and a canonical surjective map 
\[ \pi: \Fl(E) \rightarrow \Fl(N,E) .\]
The map $\pi$ is locally trivial in the Zariski topology with fibers 
\[ \pi^{-1}(s)\cong \Fl_1 \times \cdots \times \Fl_l, \]
where $\Fl_i$ is the complete flag variety of an $n_i$-dimensional $k$-vector space. Sometimes this fibration is used to reduce the study of the partial flag variety to the study of the complete flag variety.
There are similar constructions valid in the relative situation for flag bundles.

\begin{example}  Special values of L-functions and Beilinson's conjectures.\end{example}

Let $X$ be a quasi projective scheme of finite type over $\Z$. In \cite{nekovar}, Section 6 the notion of a \emph{regulator map}
\[  r_X: \K'_{2j-i}(X)_{\Q}^{(j)}     \rightarrow   \operatorname{H}_{MX_{\mathbb{R}} }^i(X_{\mathbb{R}},j)    \]
is defined, where $ \operatorname{H}_{MX_{\mathbb{R}} }^i(X_{\mathbb{R}},j)$ is \emph{motivic cohomology} of $X_{\mathbb{R}}:=X\otimes_{\Z} \mathbb{R}$.
In Conjecture $6.1-6.5$ in \cite{nekovar} precise conjectures are stated relating special values of $\L(X,s)$ at integers to the map $r_X$. These conjectures are referred to as the \emph{Beilinson conjectures}.

In \cite{gillet} the author defines for any cohomology theory $\H^*_{\alpha}(- ,\Z(i))$ satisfying a set of axioms, and any quasi projective scheme $X$ of finite type over $\Z$ Chern class maps
\[ c_i: \K'_m(X)_{\Q}^{(i)} \rightarrow \H^{2i-m}_{\alpha}(X, \Z(i)) .\]
When $\alpha:=D$ and $\H^*_{D}$ is Deligne cohomology we get Chern character maps
\[ ch_i: \K'_m(X)_{\Q}^{(i)} \rightarrow \H^{2i-m}_{D}(X\otimes \C, \mathbb{R}(i))^{+} .\]
The Chern character map $ch_i$ is a regulator map for Deligne-Beilinson cohomology, and the map $ch_i$ has been used by Borel in \cite{borel1} to prove the Beilinson conjectures for the ring $\O_K$ when $K$ is any number field. 
There are the well known formulas for the values of the Riemann zeta function included in any elementary course in calculus and integration:
\[ \L(\Z,2):=\sum_{n=1}^{\infty} \frac{1}{n^2}=\frac{1}{6}\pi^2 ,\]
\[ \L(\Z,4):=\sum_{n=1}^{\infty}\frac{1}{n^4}=\frac{1}{90}\pi^4\]
and
\[ \L(\Z,6):=\sum_{n=1}^{\infty} \frac{1}{n^6}=\frac{1}{945}\pi^6 .\]
 
In general there are the following results:
\begin{align}
&\label{neg} \L(\Z,1-k)=-\frac{B_k}{k} 
\end{align}
with $k>0$ an integer, and
\begin{align}
&\label{pos} \L(\Z,2m)=(-1)^{m-1}\frac{(2\pi)^{2m}B_{2m}}{2(2m)!} 
\end{align}
where $m\geq 1$ an integer. The number $B_i$ is the $i'th$ Bernoulli number. The formulas in \ref{neg} and \ref{pos} go back to Euler and Riemann (see the introductory book \cite{neukirch} Section VII.1 for more information). 

One would like to check if the Chern character map $ch_i$ can be used to calculate special values of the L-function $\L(\Fl(N,E),s)$ where $\Fl(N,E)$ is any flag bundle on $\O_K$, 
generalizing of Borels formula \ref{bor} to arithmetic flag schemes in any dimension.  The Beilinson conjectures are known for rings of integers in algebraic number fields, Dirichlet L-functions, some elliptic curves, Shimura curves and Hilbert-Blumenthal surfaces. See Section 8 in the paper \cite{nekovar} for more precise information and references. 

\begin{example}  Values of L-functions of flag bundles over $\O_K$ at integers. \end{example}

Let $E$ be a free $\Z$-module of rank $n$ and let $\Fl(N,E)$ be the flag bundle of type $N$ on $S:=\Spec(\Z)$. It follows 

\[ \L(\Fl(N,E),s)=\prod_{i=0}^n \L(S, s-d_i)^{l_i}=\prod_{i=0}^n \L(\Z, s-d_i)^{l_i}   .\]
Hence
\[ \L(\Fl(N,E), 1-k)=\prod_{i=0}^n \L(\Z, 1-k-d_i)^{l_i} =(-1)^{n+1}\prod_{i=0}^n (\frac{B_{k+d_i}}{k+d_i})^{l_i} \]
for $k>0$ a positive integer.

The values of $\L(\Fl(N,E),s)$ at positive integers is by Lemma \ref{positive} determined by the values of $\L(\Z,s)$ at positive integers. 

If $K$ is a number field with ring of integers $\O_K$ and $E$ a locally trivial rank $n$ $\O_S$-module where $S:=\Spec(\O_K)$ it follows again by Lemma \ref{positive} the values of $\L(\Fl(N,E),s)$ is completely determined by the values of $\L(\O_K,s)$. 

Borel discovered in \cite{borel1} regulator maps
\[  r: \K'_{4k-1}(\Z)\otimes_{\Z} \mathbb{R} \cong \mathbb{R} .\]
A non-zero element $a \in \K'_{4k-1}(\Z)\otimes_{\Z} \Q$ maps to a well defined element $R_{2k}:=r(a)\in \mathbb{R}^*/\Q^*$. This gives a formula
\begin{align}
&\label{bor} \L(\Z,2k-1)\equiv R_{2k} \text{ mod }\Q^* .
\end{align}
Formulas similar to  \ref{bor} exist for any algebraic number field $K$ and its ring of integers $\O_K$. Borel's formula for the special values of $\L(\O_K,s)$ is defined up to multiplication with a non-zero rational number. Hence if we view the values in \ref{neg} and \ref{pos} as elements in $\mathbb{R}^*/\Q^*$, it follows \ref{neg} and \ref{pos} are recovered by the formula from \cite{borel1}. 

Recent work of Bloch and Kato give an explicit formula with values in the real numbers. Lemma \ref{positive} and Borel's formula gives an explicit formula for an element 
\begin{align}
&\label{flq} \L(\Fl(N,E),m) := \prod_{i=0}^n \L(\O_K, m-d_i)^{l_i}\in \mathbb{R}^*/\Q^* 
\end{align}
for any partial flag bundle $\Fl(N,E)$ on $\O_K$.
In \cite{bloch} the authors conjecture a formula for an element $\L(M,m)\in \mathbb{R}^*$ where $M$ is a "motive", generalizing the formula in \ref{flq}. The formula conjectured in \cite{bloch}
is known to hold for some number fields and elliptic curves with complex multiplication. By Lemma \ref{positive} it follows the study of the Bloch-Kato conjecture for partial flag bundles
is reduced to the study of rings of integers in number fields.

\section{Appendix A: The weight space decomposition for algebraic K-theory of projective bundles}

In this section we calculate explicitly the weight spaces $\K'_m(\P(E^*))_{(i)}$ for any $\P^d$-bundle on $S$ to illustrate that it is easy to make explicit calculations for projective bundles.
The calculation is not neccessary for the main results of the paper, but it shows how to perform such calculations using elementary methods. We get an explicit formula for the Euler characteristic
$\chi(\P(E^*),i)$ of any projective bundle $\P(E^*)$ on $\Spec(\O_K)$ for any number field and any integer $i$.


Let in the following $X:=\Proj(\Z[x_0,..,x_n])$ be projective n-space over the ring of integers $\Z$. By the projective bundle formula for algebraic K-theory we get

\[\K_m(X)_{\Q} = \K_m(\Z)_{\Q} \otimes \Q[t]/(t^{n+1})= \K_m(\Z)_{\Q} \otimes \Q\{1,t,..,t^n\},\]

 where $t=1-L=1-[\O(-1)]$ with $L=[\O(-1)]$ and $\O(-1)$ is the tautological bundle on projective space $X:=\P(V)$. Let $R:=\Q[t]/(t^{n+1})=\Q\{1,t,t^2,..,t^n\}$.
Let

\begin{align}
&\label{log} x:=  \ln(1-t)=-(t + (1/2)t^2 + (1/3)t^3 + … + (1/n)t^n) 
\end{align}

in the ring $R=\Q\{1,t,t^2,..,t^n\}$.

\begin{lemma} Let $\psi^k$ be the kth Adams operator acting on $R$. The following holds for all integers $k \geq 0$:

\begin{align}
&\label{ad1}\psi^k(x) = kx.\\
&\label{ad2}  \text{For every integer }i \geq 1 \text{ we get } \psi^k(x^i)=k^ix^i.
\end{align}
\end{lemma}
\begin{proof} By definition $L=[\O(-1)]$ is the class in $\K_0(X)$ of the tautological line bundle $\O(-1)$ on projective space, hence the Adams operator $\psi^k$ acts as follows:
$\psi^k(L)=L^k$. We get since $t=1-L$ the following calculation:

\[\psi^k(x)=\psi^k(-(t+(1/2)t^2+(1/3)t^3+…+(1/n)t^n)= \]
\[ \psi^k(-((1-L)+(1/2)(1-L)^2+(1/3)(1-L)^3+….+(1/n)(1-L)^n)= \]
\[ -((1-\psi^k(L))+(1/2)(1-\psi^k(L))^2+…+(1/n)(1-\psi^k(L))^n)= \]
\[ \ln(\psi^k(L))=\ln(L^k)=k\ln(L)=kx \]

by Corollary A2 in the Appendix. Claim 1 is proved. Claim 2: We get
$\psi^k(x^i)=\psi^k(x)^i=(kx)^i=k^ix^i$ and Claim 2 is proved. 
\end{proof}

Note: Formal properties of exponential power series and logarithm power series valid in the formal power series ring $\Q[[t]]$ implies similar properties for exponentials and logarithms in the quotient ring $R:= \Q[[t]]/(t^{n+1})$. Formula \ref{log} was communicated to me by Charles Weibel.

If we define

\begin{align}
&\label{logarithm} \ln(L):=\ln(1-t)=-(t+(1/2)t^2+(1/3)t^3+…+(1/i)t^i+…. )  ,
\end{align}

where $\ln(L)$ lives in the formal power series ring $\Q[[t]]$,  one proves there is an equality of formal power series $\ln(L^k)=k\ln(L)$ for all integers $k \geq 0$ in $\Q[[t]]$. For a proof of this property see the Appendix. It follows the vector $x^i$ is an eigen vector for $\psi^k$ with eigen value $k^i$.  It follows the inclusion of vector spaces

\begin{align}
 &\label{iso1} \Q\{1,x,x^2,..,x^n\} \subseteq \Q\{1,t,t^2,..,t^n\} 
\end{align}

Is an isomorphism of vector spaces: The vectors $\{1,x,x^2,..,x^n\}$ are linearly independent over $\Q$ since they have different eigenvalues with respect to $\psi^k$ - the k’th Adams operator. Hence \ref{iso1} gives a decomposition of $R:=\Q[t]/(t^{n+1})$ into eigen spaces for the Adams operations $\psi^k$ for $k\geq 0$. We get an isomorphism of abelian groups

\begin{align}
\K_*(X)_{\Q} \cong \K_*(\Z)_{\Q}\otimes_{\Q} \Q\{1,x,x^2,..,x^n\}.
\end{align}

We get the following formula for $\K_m(X)_{\Q}$:

\begin{align}
&\K_m(X)_{\Q}=0\text{  if $m<0.$}\\
&\K_m(X)_{\Q}=\Q\{1,x,x^2,..,x^n\}\text{ if $m=0.$}\\
&\K_m(X)_{\Q}=0\text{ if $m=1$ or $m=2k, k \geq 1.$}\\
&\K_m(X)_{\Q}=\Q\{1,x,x^2,..,x^n\}\text{ if $m=4k+1, k >0.$}\\
&\K_m(X)_{\Q}=0\text{ if $m=4k+3, k \geq 0.$}
\end{align}

For the field of rational numbers $\Q$ we have $r_1=1$ and $r_2=0$.

\begin{lemma} The following holds for $\K_m(X)_{\Q}^{(i)}$ and $i=0,..,n$:
\begin{align}
&\K_m(X)_{\Q}^{(i)}= 0\text{ if $m<0.$}\\
&\K_m(X)_{\Q}^{(i)}=\Q\text{ if $m=0.$}\\
&\K_m(X)_{\Q}^{(i)}=0\text{ if $m=1$ or $m=2k$ with $k \geq 1.$}\\
&\K_m(X)_{\Q}^{(i)}=\Q\text{ if $m=4k+1$ with $k>0.$}\\
&\K_m(X)_{\Q}^{(i)}=0\text{ if $m=4k+3$ with $k \geq 0.$}
\end{align}
\end{lemma}
\begin{proof} The Lemma follows from the discussion above: The basis $\{1,x,x^2,..,x^n\}$ gives a decomposition of $R:=\Q[t]/(t^{n+1)}$ into eigen spaces for the Adams operation $\psi^k$ and the Lemma follows from the projective bundle formula and the calculation of $\K_m(\Z)_{\Q}$ given above. 
\end{proof}

\begin{corollary} For all $m=4k+1$ with $k>0$ and all $i=0,..,n$ it follows $\K_m(X)_{\Q}^{(i)} = \Q \neq 0$.
\end{corollary}
\begin{proof} This follows from Lemma 1 above. \end{proof}

Algebraic K-theory $\K_m(\O_K)_{\Q}$ is well known, the Adams eigen space $\K_m(\O_K)_{\Q}^{(i)}$ is well known by \cite{GGKKW}, Volume 1, Theorem 47 
and the projective bundle formula holds for $\P(E^*)$:

\[ \K_*(\P(E^*))_{\Q} \cong \K_*(\O_K)_{\Q} \otimes \Q[t]/(t^{n+1}).\]

Hence the study of the eigen space $\K_m(\P(E^*))_{\Q}^{(i)}$ is by the above calculation reduced to the study of $\K_m(\O_K)_{\Q}^{(i)}$ which is well known by Theorem \ref{boreltheorem}. 
We get the following Theorem:

\begin{theorem}\label{mainthm} Let $\Q \subseteq K$ be an algebraic number field with ring of integers $\O_K$. Let $r_1,r_2$ be the real and complex places of $K$. Let $\P(E^*)$ be a rank $e$ projective bundle on
$S:=\Spec(\O_K)$ and let $\K_m(\P(E^*))_{\Q}$ denote the m’th algebraic K-theory of the category of algebraic vector bundles on $\P(E^*)$ with rational coefficients. The following holds: Let $j \geq 0$ be an integer.
\begin{align}
&\K_m(\P(E^*))_{\Q}^{(j)} = 0 \text{ for all $m<0$ and $m=2k$ with $k \geq 1$ an integer}.\\
&\K_0(\P(E^*))_{\Q}^{(j)} = \Q \text{ if }j=0,1,2,..,e.\\
&\K_0(\P(E^*))_{\Q}^{(j)} = 0 \text{ if }j > e.\\
&\K_{4a+1}(\P(E^*))_{\Q}^{(j)} = \Q^{r_1+r_2}\text{ if $j$ is in }I:={2a+1,2a+2,..,2a+1+e}.\\
&\K_{4a+1}(\P(E^*))_{\Q}^{(j)} = 0 \text{ if $j$ is not in $I$.}\\
&\K_{4a+3}(\P(E^*))_{\Q}^{(j)} = \Q^{r_2} \text{ if $j$ is in $J:={2a+2,2a+3,..,2a+2+e}$.}\\
&\K_{4a+3}(\P(E^*))_{\Q}^{(j)} = 0 \text{ if $j$ is not in $J$.}
\end{align}
Here $a \geq 0 $ is an integer.
\end{theorem}

\begin{proof} This follows from the calculation of $\K_m(\O_K)_{\Q}^{(j)}$, the projective bundle formula and the eigen space decomposition
$R:= \Q[t]/(t^{e+1}) = \Q\{1,x,x^2,..,x^e\}$
of the ring $R$, with $x:=\ln(L):=\ln(1-t) \in R$, as described above. 
\end{proof}

\begin{corollary} \label{eulerequal} Let $S:=\Spec(\O_K)$ with $K$ a number field and let $E$ be a locally trivial $\O_S$-module of rank $d+1$. Let $F:=\O_S^{d+1}$ be the free $\O_S$-module of rank $d+1$.
The following holds:
\begin{align}
\label{proj1}&\K_m'(\P(E^*))_{\Q}^{(i)}=\K_m'(\P(\O_S^{d+1}))_{\Q}^{(i)}\text{ for any integer $i$.}\\
\label{proj2}&\chi(\P(E^*),i)=\chi(\P(\O_S^{d+1}),i)\text{ for any integer $i$.}
\end{align}
\end{corollary}
\begin{proof} By the projective bundle formula for algebraic K-theory it follows
\[ \K_m'(\P(E^*))_{\Q}\cong \K_m'(\O_S)_{\Q}\otimes \Q[t]/(t^{d+1}) ,\]
and the element $x^i$ constructed above is an eigenvector for $\psi^k$ with eigenvalue $k^i$. We get the following calculation:
\[ \K_m'(\P(E^*))_{\Q}^{(i)}=\oplus_{u+v=i} \K_m'(\O_S)_{\Q}^{(u)}x^v=\K_m'(\P(\O_S^{d+1}))_{\Q}^{(i)} .\]
Hence equation \ref{proj1} is proved. We get for any integer $k$ the following holds:
\[ \chi(\P(E^*),i)= \sum_{m\geq 0}(-1)^{m+1}dim_{\Q}(\K_m'(\P(E^*))_{\Q}^{(i)}) =\]
\[ \sum_{m \geq 0}(-1)^{m+1}dim_{\Q}(\K_m'(\P(\O_S^{d+1})_{\Q}^{(i)})=\chi(\P(\O_S^{d+1}),i) \]
and equation \ref{proj2} holds. The Corollary is proved.
\end{proof}

Note: Corollary \ref{eulerequal} gives an elementary and explicit proof of the Soule conjecture for projective space bundles on $\Spec(\O_K)$ (see Theorem \ref{lfunctpbundle}).

\begin{example} Example of Theorem \ref{mainthm} for terms $m=0,1,2,3$. \end{example}

$m=0:$
\[\K_0(\P(E^*))_{\Q}^{(l)} = \Q \text{ if $l=0,1,2,..,e.$}\]
\[\K_0(\P(E^*))_{\Q}^{(l)} = \Q \text{ if $l > e.$}\]
$m=1$:
\[\K_1(\P(E^*))_{\Q}^{(l)} = \Q^{r_1+r_2-1} \text{ if $l=1,2,3,..,e+1.$}\]
\[\K_0(\P(E^*))_{\Q}^{(l)} = 0 \text{ if $l=0$ or  $l>e+1.$}\]
$m=2$:
\[\K_2(\P(E^*))_{\Q}^{(l)}=0.\]
$m=3$:
\[\K_3(\P(^*))_{\Q}^{(l)} = \Q^{r_2} \text{ if $l=2,3,4,..,e+2.$}\]
\[\K_3(\P(E^*))_{\Q}^{(l)} = 0 \text{ if $I \neq 0,1$ or $l > e+2.$}\]

\begin{example} Schubert calculus for algebraic K-theory.\end{example}

 In a future paper a similar theory and calculation will be developed for the algebraic K-theory $\K_*(\mathbb{G}(m,E))$ of the grassmannian $\mathbb{G}(m,E)$ of $E$.
The aim of this study is to introduce and study \emph{Schubert calculus} for the K-theory of the grassmannian and flag schemes $\mathbb{F}(E)$ of a bundle $E$ over $S:=\Spec(\O_K)$, and to relate this study 
to Bloch's higher Chow groups. In \cite{sga6}, Proposition 3.1 (Berthelot's talk) the following formula is proved: Let $S$ be a noetherian scheme, $E$ a locally trivial $\O_S$-module of rank $n$ and
$P:=(p_1,..,p_k)$ a set of positive integers with $\sum_i p_i=n$ and $\mathbb{F}_P(E):=\mathbb{F}(P,E)$ the flag bundle of $E$ of type $P$, it follows the canonical morphism
\[ \K_*(S)\otimes_{\K^*(S)} \K^*(\mathbb{F}_P(E)) \cong \K_*(\mathbb{F}_P(E)) \]
is an isomorphism. Hence a formula similar to the projective bundle formula is known for flag bundles. One wants to calculate weight space decomposition
\[ \K_m(\mathbb{F}_P(E))_{\Q}  \cong \oplus_{i\in \Z}\K_m(\mathbb{F}_P(E))_{(i)} \]
for all integers $m$.

\begin{corollary} Let $X$ be a scheme of finite type over $\Spec(\O_K)$. There are no integers $M,L >>0$ with the property that $\K_m(X)_{\Q}^{(l)}=0$ for $m \geq M$ and $l \geq L$.
\end{corollary}
\begin{proof} Choose an integer $a$ such that $ 2a+1 \geq max\{M,L\}$. It follows from Theorem \ref{mainthm} that
$\K_{4a+1}(\P(E^*))_{\Q}^{(2a+1)}={\Q}^{r_1+r_2} \neq 0$. By choice $4a+1 \geq M$ and $2a+1 \geq L$.
\end{proof}


\section{Appendix B: Some general properties of formal power series}

In this section we recall some well known elementary facts on formal powerseries, logarithm power series and maps of abelian groups.

Recall the following results from [Bour], page A.IV.39 on formal power series: Let

\[ l(g(x)):= \sum_{n \geq 1}(-1)^{n-1}(1/n)(g(x))^n \in \Q[[x]]. \]

For any $g(x) \in \Q[[x]]$. Define the following formal power series:
\[ \Log(g(x)) := l(g(x)-1) \]

For any power series $g(x) \in \Q[[x]]$.
It follows
\[ \Log(1-x)=l(-x)= -(x+(1/2)x^2 + (1/3)x^3 + (1/4)x^4 + \cdots) \in \Q[[x]] .\]

Let $A$ be a commutative unital ring containing the field $\Q$ of rational numbers. Let $nil(A)$ be the nilradical of $A$. Let $1-nil(A)$ denote the set of elements on the form $1-u$ with $u \in nil(A)$. It follows $1-u$ is a multiplicative unit in $A$. The set $1-nil(A)$ has a multiplication: $(1-u)(1-v)=1-u-v+uv=1-z$ with $z=-u-v+uv$, and the element $z$ is again in $nil(A)$. Hence $(1-u)(1-v)=1-z$ is in $1-nil(A)$. It follows $1-nil(A)$ is a subgroup of the multiplicative group of units in $A$.

\begin{lemma} (A1) Let $u \in nil(A)$ be an element with $u^{k+1}=0$. Define the following map:
\[ \ln : 1-nil(A) \rightarrow nil(A) \]
by
\[ \ln(1-u):= -(u+(1/2)u^2 + (1/3)u^3 + \cdots + (1/k)u^k) \in nil(A).\]
It follows $\ln$ is a morphism of groups: For any two elements $1-u,1-v \in 1-nil(A)$ it follows
\[\ln((1-u)(1-v))=\ln(1-u) + \ln(1-v).\]
\end{lemma}
\begin{proof} From \cite{Bour}, page A.IV.40 we get

\[ \Log(1-x)=l(1-x-1)=l(-x)\text{ in }\Q[[x]]. \]

The   following holds for the powerseries $l(x)$:
$l(x+y+xy)=l(x)+l(y)$ in the ring $\Q[[x,y]]$. We may for any two elements $u,v$ in $nil(A)$ define a map
\[ f:\Q[[x,y]] \rightarrow A \]

 by $f(x)=u,f(y)=v$. It follows $f$ induce a well defined map of rings

\[ f': \Q[[x,y]]/I \rightarrow A\]

 where $I=ker(f)$. In the ring $\Q[[x,y]]$ we get the following formula:

\begin{align}
&\label{ln} \Log((1-x)(1-y))=\Log(1-x-y+xy):=l(-x-y+xy)= 
\end{align}
\[  l(-x-y+(-x)(-y)=l(-x)+l(-y)=\Log(1-x)+\Log(1-y). \]

It follows the same formula \ref{ln} holds in the quotient ring $\Q[[x,y]]/I$. Hence we get the following formula for the map $\ln$ (viewing u and v as elements in the quotient $\Q[[x,y]]/I$):
\[ \ln((1-u)(1-v))=\Log((1-x)(1-y))=\Log(1-x)+\Log(1-y)=\ln(1-u)+\ln(1-v). \] 
Hence the map $\ln$ is a map of groups. 
\end{proof}

Note: Lemma A1 may also be proved using Bell polynomials.

\begin{corollary} (A2) Use the notation from Lemma A1. If $1-u \in 1-nil(A)$ the following holds for any integer $k \geq 1$:
\[\ln((1-u)^k)=k\ln(1-u).\]
\end{corollary}
\begin{proof} This follows from Lemma A1 and an induction. \end{proof}

Example: Let $A:=\Q[t]/(t^{e+1})$ with $nil(A)=(t)$ define the following “logarithm” map ($u \in nil(A)$):

\[ \ln(1-u):= -(u+(1/2)u^2+(1/3)u^3 + \cdots + (1/e)u^e) \in A. \]

It follows
\begin{align}
&\label{lnp}  \ln((1-u)^k)=k\ln(1-u) 
\end{align}

 for any integer $k \geq 1$. The property \ref{lnp} is well known when we consider the logarithm function defined for real numbers, and the above section proves it holds for formal power series.

Note: Formal properties of exponentials and logarithms in $\Q[[t]]$ can also be proved using Bell polynomials.

\textbf{Acknowledgements}: Thanks to Shrawan Kumar,  Marc Levine, Chris \Soule and Charles Weibel for answering questions and providing references on algebraic K-theory and flag varieties. Thanks also to Alexander Beilinson and Christopher Deninger
for answering questions on L-functions and the Beilinson and Bloch-Kato conjectures.

\end{document}